\documentclass[final,12pt]{amsart}

\usepackage{multirow}
\usepackage[centering,text={16cm,20cm},
		marginparwidth=18mm]{geometry}

\usepackage{amsmath,bm,graphicx,amsthm}
\theoremstyle{definition}
\newtheorem{theorem}{Theorem}

\newtheorem{lemma}[theorem]{Lemma}
\newtheorem{proposition}[theorem]{Proposition}

\newtheorem{definition}[theorem]{Definition}

\newtheorem{remark}[theorem]{Remark}
\numberwithin{equation}{section}

\def\Z{\mathbb{Z}}

\def\F{\mathbb{F}}

\def\R{\textbf{R}}
\def\BB{\textbf{B}}

\def\S{\mathbf{S}}

\def\x{\mathbf{x}}
\def\y{\mathbf{y}}
\def\z{\mathbf{z}}
\def\w{\mathbf{w}}

\def\QQ {\mathcal{Q}}
\def\RR {\mathcal{R}}
\def\SS {\mathcal{S}}
\def\FF {\mathcal{F}}

\def\L {\mathcal{L}}

\def\J {\mathcal{J}}
\def\I {\mathcal{I}}

\def\X{\mathbb{X}}
\def\O{\mathbb{O}}

\def\pent{{Pent}_{\t \beta \gamma}(\x, \y)}

\def\hex{{Hex}_{\t \beta \gamma \t \beta}(\x, \y)}

\def\Ctq{{C}_Q}

\def\a{\mathfrak{a}}
\def\b{\mathbf}
\def\t{\widetilde}
\def\h{\widehat}
\def\A{\alpha}
\def\B{\beta}

\def\{{\lbrace}
\def\}{\rbrace}

\newcommand{\exact}[3]{$0 \longrightarrow #1 \longrightarrow #2 \longrightarrow #3~\longrightarrow~0 $}

\begin{document}


\title{COMBINATORIAL KNOT FLOER HOMOLOGY AND DOUBLE BRANCHED COVERS}

\author{Fatemeh Douroudian}%

\address{Department of Mathematics, Faculty of Mathematical Sciences, Tarbiat Modares University, P.O. Box 14115-137, Tehran, Iran.}%
\email{douroudian@modares.ac.ir}%



\begin{abstract}
Using a Heegaard diagram for the pullback of a knot $K \subset S^3$ in its double branched cover $\Sigma_2(K)$, we give a combinatorial proof for the invariance of the associated knot Floer homology over $\Z$.
\end{abstract}

\maketitle


\section{Introduction} \label{intro}

Heegaard Floer homology, introduced by Ozsv\'ath and Szab\'o, is a collection of invariants for closed oriented three-manifolds. Various versions of Heegaard Floer homology is defined by counting some holomorphic disks in the symmetric product of a Riemann surface. There is a relative version of the theory for a pair $(Y, K)$, where $K$ is a nullhomologous knot in the three-manifold $Y$, which is developed by Ozsv\'ath and Szab\'o \cite{OS} and independently by Rasmussen \cite{Ras}. However working with holomorphic disks makes it difficult to compute the homolgy groups. In~\cite{Sarkar} Sarkar and Wang gave an algorithm that makes the computation of hat version combinatorial. In~\cite{nice}, Ozsv\'ath, Stipsicz and Szab\'o gave a combinatorial algorithm for constructing $\h{HF}(Y)$ and provide a topological proof of its invariance.

Given a knot $K\subset S^3$, a \emph{grid diagram} $G$ associated with $K$ is an $n\times n$ planar grid, together with two sets $\X=\{ X_i\} _{i=1}^{i=n}$ and $\O=\{ O_i\} _{i=1}^{i=n}$ of basepoints. Each column and each row contains exactly one $X$ and one $O$ inside. We view this grid diagram as a torus $T^2 \subset S^3$ by standard edge identifications. Here each horizontal line is an $\alpha$ circle and each vertical line is a $\beta$ circle, and $(T^2, \b \A, \b \B, \O, \X)$ is a multi-pointed Heegaard diagram for $(S^3, K)$. Manolescu, Ozsv\'ath and Sarkar~\cite{MOS} showed that such diagrams can be used to compute $\widehat{HFK}(S^3, K)$ combinatorially. In \cite{MOST}, Manolescu, Ozsv\'ath, Szab\'o and Thurston gave a combinatorial proof of the invariance of knot Floer homology and a combinatorial presentation of the basic properties of link Floer homology over $\Z$.
In \cite{L}, Levine gave a construction of a Heegaard diagram for the pullback $\t K$ of a knot $K \subset S^3$ in its $m$-fold cyclic branched cover $\Sigma_m(K)$ to compute $\widehat{HFK}(\Sigma_m(K),\t K)$ over $\Z_2$ combinatorially.
In this paper we use a recent work of ~Ozsv\'ath, Stipsicz and Szab\'o \cite{Sign}, where they assign signs to the rectangles and bigons in a nice Heegaard diagram, and we provide a combinatorial proof of the invariance of the knot Floer homology of $(\Sigma_2(K),\t K)$ with $\Z$ coefficients. Note that the invariance of $\h{HFK}(\Sigma_2(K),\t K;\Z)$ was known by the works of \cite{OS} and \cite{Ras}, but in this paper we give a combinatorial argument which allows us to compute $\h{HFK}(\Sigma_2(K),\t K;\Z)$ combinatorially.

\begin{figure}[h]
\centerline{\includegraphics[scale=0.8]{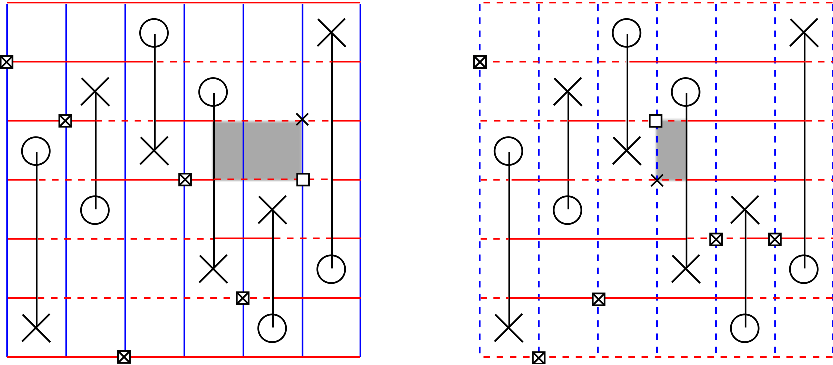}}
\caption {A Heegaard diagram $\t G=(\t T, \bm{\t \alpha}, \bm{\t \beta}, \O, \X)$ for $\t K \subset \Sigma_2(K)$, where $K$ is the figure eight knot. Here the horizontal lines represent arcs in the lifts of $\A$-circles and each of vertical lines represent a lift of a $\B$-circle. A rectangle in $\RR(\t G)$ is illustrated which contributes in the differential of $C(\t G)$ and connects a generator $\x$ which is shown with crosses to a generator $\y$ which is shown with hollow squares.}
\label{grid}
\end{figure}

The construction of a nice Heegaard diagram for $(\Sigma_m (K),\t K)$, described in \cite{L} is as follows: Let $T^2$ be the torus which describes the grid diagram $G$ of $K$. First isotope $K$ so that it lies in $H_{\A}$, the handlebody corresponding to $\A$ curves. Next let $F$ be a Seifert surface that is contained entirely in a ball in $H_{\A}$. We can see, after isotoping back both of $K$ and $F$ to $H_{\B}$, that the intersections of $F$ with $T^2$ are the vertical lines that connect the $X$ and the $O$ in each column in $G$. We define $\t T$ to be the surface obtained by gluing together $m$ copies of $T^2$, denoted by $T_0, \dots, T_{m-1}$, along the branch cuts connecting the pairs of $X$ and $O$ in each column such that whenever $X$ is above $O$, the left side of the branch cut in $T_k$ is glued to the right side of the same cut in $T_{k+1}$; if the $O$ is above the $X$, then glue the left side of the branch cut in $T_k$ to the right side of the same cut in $T_{k-1}$ (indices modulo $m$). The projection map $\pi :\t T \rightarrow T^2$ is an $m$-fold cyclic branched cover, branched over the basepoints. Each $\alpha$ and $\beta$ circle in $T^2$ bounds a disk away from branch points in $S^3-K$ so each of them has $m$ distinct lifts to $\Sigma_m(K)$, the $m$-fold cyclic branched cover. Each lift of each $\alpha$ circle intersects exactly one lift of each $\beta$ circle. This can be shown as in Fig.~\ref{grid} with $m$ disjoint \emph{grids}. Denote by $\t \B _{j}^{i}$ for $i=0, \dots, m-1$ and $j=0, \dots, n-1$ the vertical arcs, and denote by $\t \A _{j}^{i}$ the arc which has intersection with $\t \B _{0}^{i}$. Let us denote this Heegaard diagram by $\t G=(\t T,\t{\bm{\A}},\t{\bm{\B}},\O, \X)$, where $\t{\bm{\A}}$ and $\t{\bm{\B}}$ are the lifts of ${\bm{\A}}$ and ${\bm{\B}}$ to $\t T$.
In the figures for $m=2$, we draw the lifts with superscript zero (i.e. $\t\A^0_i$ and $\t\B^0_j$) with solid lines, and the lifts with superscript one (i.e. $\t\A^1_i$ and $\t\B^1_j$) with dashed lines.

Denote by $\RR(G)$ the set of embedded rectangles in $T^2$ which do not contain any basepoints and whose left and right edges are arcs of $\B$ circles and whose upper and lower edges are arcs of $\A$ circles. Each rectangle in $\RR(G)$ has $m$ disjoint lifts to $\t T$ (possibly passing through the branch cuts), denote the set of such lifts by $\RR(\t G)$.

Let the set of generators $\S(\t G)$ be the set of unordered $mn$-tuples $\x$ of intersection points between $\t \alpha^i_j$ and $\t \beta^i_j$ for $i = 0 , \dots, m-1$ and $j= 0, \dots, n-1 $, such that each of $\t \alpha^i_j$ and $\t \beta^i_j$ has exactly one component of $\x$.  Also denote by $\S(G)$ the set of generators for $G$. In \cite{L}, Levine showed that any generator $\x \in \S(\t G)$ can be decomposed (non-uniquely) as $\x= \t\x_1\cup \dots \cup \t\x_m$, where $\x_1, \dots, \x_m$ are generators in $\S(G)$, and $\t\x_i$ is a lift of $\x_i$ to $\t G$.

The Alexander grading of a generator $\x \in \S(\t G)$ is defined as follows: Given two finite sets of points $A$, $B$ in the plane, let $\I (A, B)$ be the number of pairs $(a_1,a_2)\in A$ and $(b_1,b_2)\in B$ such that $a_1<b_1$ and $a_2<b_2$. Let $\J (A, B)=\dfrac{1}{2}(\I (A, B)+\I (B,A))$. Given $\x_i \in \S(G)$, define 
\begin{center}
$A(\x_i)= \J(\x_i-\dfrac{1}{2}(\X+\O),\X-\O)- (\dfrac{n-1}{2})$.
\end{center}
in which we consider $\J$ as a bilinear function of its two variables. For any generator $\x \in \S(\t G)$, consider one of the decompositions $\x= \t\x_1\cup \dots \cup \t\x_m$ and define
\begin{center}
$A(\x)= \dfrac{1}{m} \sum_i A(\x_i)$.
\end{center}
A simple calculation shows that the value of $A(\x)$ in the above definition is well-defined (i.e. it is independent of the choice of the decomposition).  

Let $C({\t G})$ be the $\Z_2$ vector space generated by $\S(\t G)$. Define a differential $\partial$ on $C({\t G})$ by setting a nonzero coefficient for $\y$ in $\partial \x$ if and only if $\x, \y \in \S(\t G)$ agree along all but two vertical circles and there is a rectangle $R \in \RR(\t G)$ whose lower-left and upper-right corners are in $\x$ and whose lower-right and upper-left corners are in $\y$, and which does not contain any components of $\x$ in its interior. We denote by $Rect(\x, \y)$ the set of such rectangles.

\begin{center}
$\partial \x=\sum\limits_{\y \in \S(\t G)} \sum\limits_{R\in Rect(\x,\y)} \y$
\end{center}

Note that the $\partial$ preserves the Alexander grading. 

Let $H_*(\t G)$ be the homology of the chain complex $(C(\t G), \partial)$.

Later we consider $C(\t G)$ as a free $\Z$-module and define the boundary map using the signs that are assigned to each $R\in Rect(\x,\y)$.

\begin{remark}\label{horver}
In this construction we cut along the vertical lines connecting the $X$ and the $O$ in each column, and glue two copies of the resulting surface.
Note that we can alternatively consider a different Seifert surface, with the property that its intersection with the torus $T^2$ be the horizontal lines connecting the $X$ and the $O$ in each row, and glue two copies of the resulting surface.
In this way we get the same Heegaard diagram up to homeomorphism (this follows for example from \cite{Lickorish} Theorem 7.9), and so in the definition of $(C(\t G), \partial)$ we could have used this alternative construction and it would not change the resulting complex, and in priori the homology groups.
\end{remark}

Other useful notations are as follows: Let $\x , \y \in \S(\t G)$. A \emph{path} $\gamma$ from $\x$ to $\y$ is a closed oriented path composed of arcs on $\A -$ and $\B -$circles, which has its corners among $\x$ and $\y$, oriented so that $\partial(\gamma \cap \A)=\y -\x$ and $\partial(\gamma \cap \B)=\x -\y$. Let $D_1 , \dots , D_m$ be the components of $\t T \backslash {\bm \A \cup \bm \B}$, a \emph{domain} $p$ from $\x$ to $\y$ is a two-chain $D = \sum a_i D_i$ in $\t T$ whose boundary $\partial p$ is a path from $\x$ to $\y$. Define $\pi (\x, \y)$ to be the set of domains from $\x$ to $\y$. There is a natural composition law 

\begin{center}
$* : \pi(\b a,\b b) \times \pi(\b b,\b c) \longrightarrow \pi(\b a,\b c)$
\end{center}


In order to state the main theorem of this paper, we need one more definition. Following \cite{OSS}, we define the stable knot Floer homology $\h{HFK}_{st}(Y,K)$ of a knot inside a three manifold $Y$ as follows:

\begin{definition}
Two pairs $(V_1, a_1)$ and $(V_2, a_2)$ of vector spaces, where $V_1$, $V_2$ are finite dimensional free modules over the ring $\F$ and $a_1,a_2$ are non-negative integers, are called equivalent pairs if either $a_1\geq a_2$ and $V_1\cong V_2 \otimes (\F \oplus \F)^{a_1 -a_2}$, or $a_2\geq a_1$ and $V_2\cong V_1 \otimes (\F \oplus \F)^{a_2 -a_1}$. If we are working with graded vector spaces then we require the isomorphisms to preserve the grading.
\end{definition}

Let $D=(\Sigma, \bm{\alpha}, \bm{\beta}, \w, \z)$ be a nice Heegaard decomposition of $Y$ compatible with the knot $K \subset Y$ such that $|\w |=|\z |$. Then the equivalence class of the pair $(H_*(\t G), |\z |)$ is the stable knot Floer homology of $K \subset Y$ denoted by $\h{HFK}_{st}(Y,K)$ (as graded vector space).

Note that $\h{HFK}_{st}(Y,K)$ and $\h{HFK}(Y,K)$ are in fact equivalent. As shown in \cite{OSS, L, MOST}, we have $H_*(\t G) \cong \h{HFK}(Y,K) \otimes (\F \oplus \F)^{|\z | -1} $, where $\F$ is either $\Z_2$ or $\Z$ depending on the way we define the chain complex $(C(\t G), \partial)$. The statement of the main theorem is as follows:

\begin{theorem} \label{main}
Let $K$ be an oriented knot and $\t K$ be its pullback in the double branched cover $\Sigma_2(K)$. Then the stable knot Floer homology of $\t K \subset \Sigma_2(K)$ over $\Z$ is an invariant of the knot.
\end{theorem}

In Section~\ref{assign}, we review some definitions and theorems from~\cite{Sign}. In Section~\ref{invar} we prove the invariance of the knot Floer homology of the pullback of a knot $K \subset S^3$ in its double branched cover $\Sigma_2(K)$ over $\Z$.


\section{Prelimineries} \label{assign}

In~\cite{Sign}, it is shown that there exists a sign assignment for nice Heegaard diagrams, using this fact, we want to extend the definition of boundary map from the case of $\Z_2$ to the case with coefficients in $\Z$. 
To be self-contained, we review some definitions and theorems, for more details see~\cite{Sign}.

We call an empty rectangle or bigon in a nice Heegaard diagram, a \emph{flow}. Fix two sets $\bm \A$ and $\bm \B$ such that $|\bm \A|=|\bm \B|=m$. We define formal generators and formal flows as follows:

\begin{definition}
A formal generator of power $m$ is a pair $(\epsilon,\sigma)$, where $\epsilon = (\epsilon_1, \dots , \epsilon_m) \in {\{ \pm 1 \} }^m$ and $\sigma \in S_m$. A formal generator of power $m$ has the following equivalent definition: a pairing of elements of $\bm \A$ and $\bm \B$ with an assignment of $ \pm 1 $ to each pair.
\end{definition}

We can represent this definition pictorially. Draw $m$ oriented crossings and label the oriented arcs at the $i^{th}$ crossing with $\A_i$ and $\B_{\sigma(i)}$. The sign $\epsilon_i$ is the intersection number $\A_i\cdot \B_{\sigma(i)}$, computed with respect to the standard orientation of the plane.

\begin{definition}
Consider two balls $U$ and $V$ of radii $1$ embedded in the complex plane, such that $U$ is centered at $0$ and $V$ is centered at the point $\sqrt{2}$. Let $\BB = U\cap V$, we call the point $\frac{\sqrt{2}}{2} - \frac{\sqrt{2}}{2} i \in \partial U \cap \partial V$ the initial point of $\BB$ and the point  $\frac{\sqrt{2}}{2} + \frac{\sqrt{2}}{2} i \in \partial U \cap \partial V$ the terminal point of $\BB$. The arc $\partial U \cap \BB$ is called $\bm a$ and the arc $\partial V \cap \BB$ is called $\bm b$. We call this setting the \emph{ideal bigon}, and when there is no confusion denote it only by $\BB$.
\end{definition}

\begin{definition}
Let $\x = (\epsilon , \sigma)$ and $\y = (\epsilon ' , \sigma ')$ be two formal generators of power $m$ such that for some $i \in \{ 1, \dots , m \}$ we have $\epsilon_i = -\epsilon_i '$ and $\epsilon_j = \epsilon_j '$ for all $j\neq i$, and also $\sigma = \sigma '$. A \emph{formal bigon} of power $m$, from $\x$ to $\y$ is the data of an orientation of the arcs $\bm a$ and $\bm b$ of the ideal bigon $\BB$ with the following (equivalent) properties:
\begin{itemize}
\item the local intersection number of $\bm a$ and $\bm b$ at the initial point of $\BB$ (with respect to the orientations that we fixed) is $\epsilon_i$
\item the local intersection number of $\bm a$ and $\bm b$ at the terminal point of $\BB$ (with respect to the orientations that we fixed) is $\epsilon_i '$
\end{itemize}

We denote such formal bigon by $\mathcal{B} :\x \rightarrow \y$.
\end{definition}

\begin{definition}
Consider a square $\R$ embedded in complex plane with corners $0 , 1 , i , 1+i$. We label the upper (resp. lower) edges $N$ (resp. $S$). Also label the right (resp. left) edges $E$ (resp. $W$). We call the pair $(SW , NE)$ (of intersection poits) the initial pair of points, and $(SE , NW)$ the terminal pair of points. We call $\R$ an \emph{ideal rectangle}.
\end{definition}

\begin{definition}\label{formalrec}
Let $\x = (\epsilon , \sigma)$ and $\y = (\epsilon ' , \sigma ')$ be two formal generators of power $m$ such that for disjoint $s,n \in \{ 1, \dots , m \} $ and disjoint $e,w \in \{ 1, \dots , m \} $ we have:
\begin{itemize}
\item $\sigma(s)=w$ , $\sigma(n)=e$
\item $\sigma'(s)=e$ , $\sigma '(n)=w$
\item $\sigma(j)=\sigma '(j)$ for all $j\neq s, n$
\item $\epsilon_j = \epsilon_j '$ for all $j\neq n, s.$ 
\end{itemize}
A \emph{formal rectangle} of power $m$, from $\x$ to $\y$ which we denote by $\mathcal{R} :\x \rightarrow \y$, is the data of an ideal rectangle $\R$ with a choice of orientation for the four edges of $\R$ with following properties:
the local intersection signs for the initial points $SW$ and $NE$ coincide with the signs $\epsilon_s$ and $\epsilon_n$. Also the local intersection signs for the terminal points $SE$ and $NW$ coincide with $\epsilon'_s$ and $\epsilon'_n$. 
\end{definition}

For such a formal rectangle, we denote the edges $N$, $S$, $E$ and $W$ by $\A_n$, $\A_s$, $\B_e$ and $\B_w$. We represent such formal rectangle by the sequence $\A_n \rightarrow \B_w \rightarrow \A_s \rightarrow \B_e$. The data of $\sigma(j)$ and $\epsilon(j)$ for $j \neq n,s$ are not explicitly mentioned in this \emph{arrow notation} and should be understood from the context.

We call either a formal bigon or a formal rectangle, a \emph{formal flow}. Given two formal generators $\x = (\epsilon , \sigma)$ and $\y = (\epsilon ' , \sigma ')$, we denote by $\phi :\x \rightarrow \y$ a formal flow from $\x$ to $\y$. Denote by $\FF_m$ the set of formal flows of power $m$ (with respect to the initial and terminal formal generators). Given two formal flows $\phi_1 :\x \rightarrow \y$ and $\phi_2 :\y \rightarrow \z$, for ease of language we call the pair $(\phi_1 , \phi_2)$, the \emph{composition} of $\phi_1$ and $\phi_2$. Note that although we call the pair $(\phi_1 , \phi_2)$ the composition of $\phi_1$ and $\phi_2$, , $(\phi_1 , \phi_2)$ is not a flow and it is not of the same type as $\phi_1$ and $\phi_2$.

\begin{figure}[h]
\centerline{\includegraphics[scale=0.5]{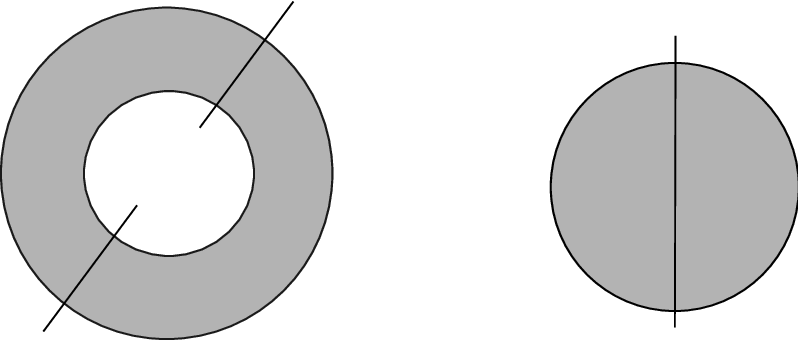}}
\caption {Boundary degenerations of composition of flows.}
\label{dege}
\end{figure}

\begin{definition}
Let $\x = (\epsilon , \sigma)$, $\y = (\epsilon ' , \sigma ')$ and $\z = (\epsilon '', \sigma '')$ be formal generators of the same power, and $\phi_1 :\x \rightarrow \y$ and $\phi_2 :\y \rightarrow \z$ be two formal flows.
If with some orientations and labelling of arcs with $\A_i$ and $\B_j$
the pair $(\phi_1 , \phi_2)$
has one of the forms in Fig.~\ref{dege}, then we say the pair $(\phi_1 , \phi_2)$ is a \emph{boundary degeneration}.
Note that in this case $\z =\x $.
The boundary degeneration is of Type $\A$ (resp. $\B$), when the circle(s) in Fig.~\ref{dege} is decorated with $\A$ (resp. $\B$).
\end{definition}

We define a sign assignment of power $m$ as follows:

\begin{definition} \label{def sign} 
A sign assignment $\SS$ of power $m$ is a map $\SS :\FF_m \longrightarrow \{\pm 1\}$ that satisfies the following conditions: 

(S-1) if the composite flow $(\phi_1 , \phi_2)$ is a degeneration of Type $\A$, then 
\begin{center}
$\SS(\phi_1) \cdot \SS(\phi_2)=1$
\end{center}

(S-2) if the composite flow $(\phi_1 , \phi_2)$ is a degeneration of Type $\B$, then 
\begin{center}
$\SS(\phi_1) \cdot \SS(\phi_2)= -1$
\end{center}

(S-3) given two pairs $(\phi_1, \phi_2)$ and $(\phi_3, \phi_4)$ such that the initial formal generator of $\phi_1$ and $\phi_3$ are the same and the terminal formal generator of $\phi_2$ and $\phi_4$ are the same, then 
\begin{center}
$\SS(\phi_1) \cdot \SS(\phi_2)+\SS(\phi_3) \cdot\SS(\phi_4)=0$.
\end{center}
\end{definition}

We can construct new sign assignments from an old one:

\begin{definition}
Let $\SS$ be a sign assignment and $u$ be any map from the set of formal generators of power $m$ to $\{\pm 1\}$, then we can define a new sign assignment $\SS'$ such that for any $\phi:\x\rightarrow\y$ in $\FF_m$, $\SS'(\phi)=u(\x)\cdot\SS(\phi)\cdot u(\y)$.
If $\SS$ and $\SS'$ are related as above, we say that $\SS$ and $\SS'$ are gauge equivalent sign assignments.
\end{definition}

Using these definitions the precise statement of the result of \cite{Sign} that we need is as follows:
\begin{theorem}
For a given power $m$ there exists a sign assignment, and it is unique up to gauge equivalence.
\end{theorem}

If a nice Heegaard diagram has $m$ $\A$-curves, we say that the Heegaard diagram is of \emph{power $m$}. Let $\mathcal{D}=(\Sigma, \bm{\A}, \bm{\B}, \b w, \b z)$ be a nice Heegaard diagram of power $m$. Denote by $\S$ the set of generators. We fix an ordering for $\bm{\A}$ and $\bm{\B}$, and an orientation for each $\A$- and $\B$-circle. Then a generator $\x$ of the Heegaard diagram specifies a formal generator $\x_f$ of power $|\bm{\A}|$. Let $\x,\y \in \S$, an empty bigon or empty rectangle $\phi$ from $\x$ to $\y$ determines a formal bigon or formal rectangle $F(\phi):\x_f \rightarrow \y_f$ of power $|\bm{\A}|$.
Fix a sign assignment $\SS$ of the same power.
For each $\x\in \SS$, the boundary operator $\t{\partial}^{\Z}(\x)$ is defined as follows:
$$\t{\partial}^{\Z}(\x)=\sum_{\y\in \S}\sum_{\phi \in \mathrm{Flows}(\x,\y)} \SS(F(\phi))\y ,$$
\noindent where we denote by $\mathrm{Flows}(\x,\y)\subset \pi(\x,\y)$ the set of empty bigons and empty rectangles from $\x$ to $\y$. The following theorem is proved in \cite{Sign}.

\begin{theorem} \label{signassign}
\cite{Sign} The map $\t{\partial} ^\Z$ over $\Z$ satisfies $(\t{\partial}^\Z)^2=0$. The resulting Floer homology $\t{HF}(\mathcal{D};\Z)$ is independent of the choice of $\SS$, the order of the $\A$- and $\B$-curves, and the chosen orientation on each of the $\A$- and $\B$-curves. 
\end{theorem}

\section{Invariance of Knot Floer Homology} \label{invar}  
In this section we want to show that the combinatorial knot Floer homology of $\t K \subset \Sigma_2(K)$ is independent of $G$, the grid diagram for $K \subset S^3$. As a result of the work of Cromwell \cite{Cromwell}, any two grid diagrams of a knot $K \subset S^3$ can be connected by a sequence of three elementary moves, each resulting in a new Heegaard diagram for $\t K \subset \Sigma_2(K)$, as follows: 

\begin{enumerate}\label{moves}
\item $\mathbf{Cyclic\;Permutation}$ This move corresponds to cyclically permuting the rows or the columns of $G$, the grid diagram for the knot $K$, and obtaining a new grid diagram $H$ for $K$. Consequently we obtain a new Heegaard diagram $\t H$ from $\t G$ for $\widetilde{K}$.
    
\item $\mathbf{Commutation}$ Consider two consecutive columns in a grid diagram $G$ of a knot $K$, the $X$ and $O$ decorations of one of the columns separate the vertical circle into two arcs. If both of the $X$ and $O$ decorations of the other column is in one of the arcs, switching these two columns is a commutation move for $G$ and leads to another grid diagram $H$ for $K$. For the knot $\widetilde{K}$, we consider the Heegaard diagram $\widetilde{H}$ to be obtained from $\widetilde{G}$ by a commutation move. Commutation can be alternatively done by reversing the roles of rows and columns.      

\item $\mathbf{Stabilization/Destabilization}$
Let $G$ be a $n\times n$ grid diagram of the knot $K$. We want to add two consecutive breaks in $K$ and introduce a new $(n+1)\times (n+1)$ grid diagram $H$ for $K$. Label the decorations of $G$ by ${\lbrace O_i \rbrace }^{i=n+1}_{i=2}$ and ${\lbrace X_i \rbrace }^{i=n+1}_{i=2}$. Consider a row with decorations $O_i$ and $X_i$ in the grid diagram $G$, split the row so that it becomes two rows, and introduce a new column somewhere between $O_i$ and $X_i$.
We copy $O_i$ onto one of the new rows and the $X_i$ in the other row. Then we add $O_1$ and $X_1$ in the two squares which are the intersections of these two rows with the new column, such that $O_1$ is in the same row as $X_i$.
We let $\A_1$ denote the new horizontal circle in $H$ which separates $O_1$ from $X_1$. The number of $\B$-circles in $H$ is one more than $G$.
We denote the $\B$-circle which is just to the left of $O_1$ and $X_1$ by $\B_1$, and the $\B$-circle just to the right of $O_1$ and $X_1$ by $\B_2$.
We call $H$ a stabilization of the grid diagram $G$ of $K$. There is a similar stabilization move where the roles of the rows and the columns are interchanged. Using commutation moves, we consider only certain stabilization moves where three of $O_1$, $X_1$, $O_i$ and $X_i$ are placed such that they have a corner in common. For the knot $\widetilde{K}$, we consider the Heegaard diagram $\widetilde{H}$ to be obtained from $\widetilde{G}$ by stabilization.
The grid diagram $\t H$ has two $\B$-circles and two $\A$-circles more than $\t G$. Destabilization is the reverse of the stabilization move.

Note that the definition of the complex is invariant if we rotate everything $180^{\circ}$. This means that without loss of generality we assume that the $O_1$ is directly above $X_1$.

\end{enumerate} 
   
\begin{remark} \label{horver2}
Note that for each elementary move, we have to consider two separate cases. One for applying the move to rows, and the other for applying it to the columns. By Remark~\ref{horver} the chain complex is independent of the choice of the direction of the cuts, so by symmetry we can assume that the cuts are vertical and we always apply the elementary moves to the columns.
\end{remark}
   
\subsection{Cyclic Permutation}

The case when we cyclically permute the columns and the complex is considered with $\Z_2$ is tautological, since in this case the complex only depends on the topology of the Heegaard diagram, and the cyclic permutation does not change the position of branch cuts. In the case with coefficient in $\Z$ one should note that the cyclic permutation still does not change the topology of the Heegaard diagram, but it will permute the order of $\alpha$ and $\beta$ curves. Since the homology is independent of this ordering (see Theorem~\ref{signassign}) we see that the complex and the $\partial$ operator will be invariant under the cyclic permutation. 

Using Remark~\ref{horver2} the case that we apply the cyclic permutation to the rows follows from the case that we apply the cyclic permutation to the columns, and so we are done.


\subsection{Commutation} \label{comm}
Let $\t G$ be a Heegaard diagram for $\t K$, and let $\t H$ be the Heegaard diagram obtained by commutation.

\begin{figure}[h]
\centerline{\includegraphics[scale=0.7]{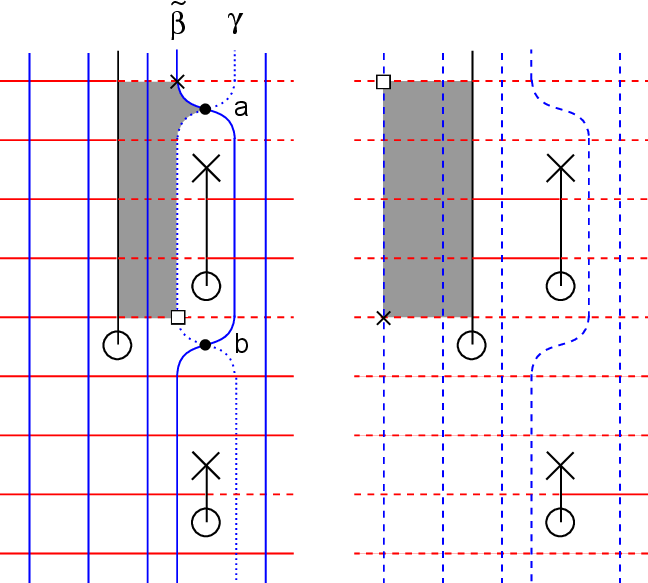}}
\caption {Here we showed two Heegaard diagrams $\t G$ and $E$ and a pentagon in $\pent$, where components of $\x$ are indicated by crosses and those of $\y$ are indicated by hollow squares.}
\label{pentagon}
\end{figure}

We define a Heegaard diagram $E$ associated to the specific commutation, consisting of two $n\times n$ grids where the opposite sides of each grid are identified. The $X$ and $O$ decorations in the right grid are the same as the decorations for $G$, in the left grid we just interchange decorations of the two columns where we want to make the commutation move. See Fig.~\ref{pentagon}. Let $\t \beta^i_j$ denote the vertical arcs in the Heegaard diagram for $i = 0, 1$ and $j= 0, \dots, n-1$. For the horizontal arcs, we denote by $\t \alpha^i_j$ the arc which has intersection with $\t \beta^i_0$.

Let the set of generators $\S(E)$ be the set of unordered $2n$-tuples $\x$ of intersection points between $\t \alpha^i_j$ and $\t \beta^i_j$ for $i= 0, 1$ and $j= 0, \dots, n-1 $, such that each of $\t \alpha^i_j$ and $\t \beta^i_j$ has exactly one component of $\x$.

Define the set of allowed regions $\RR(E)$ to be the set of rectangles that are topological disks and whose upper and lower edges are arcs in $\alpha$-circles and whose left and right edges are arcs of $\beta$-circles.

Let $C(E)$ be the free $\Z$-module generated by $\S(E)$. Given $\x, \y \in \S(E)$ such that they agree along all but two vertical circles, let $Rect_E(\x,\y)$ denote the set of rectangles $R \in \RR(E)$ whose lower-left and upper-right corners are in $\x$ and whose lower-right and upper-left corners are in $\y$, and which does not contain any $X$, $O$, or components of $\x$ in its interior. Define the differential $\partial_E$ on $C(E)$ as follows:
$$\partial_E \x=\sum_{\y\in \S(E)}\sum_{R \in Rect_E(\x,\y)} \SS(F(R)) \y ,$$
where $F(R)$ denotes the formal rectangle associated with $R$. 

It is convenient to draw the new Heegaard diagram $E$ in the same diagram as of the Heegaard diagram of $\t G$, replacing the distinguished vertical circle $\t \beta$ in $\t G$ with a different one $\gamma$ in $E$. The circles $\t \beta$ and $\gamma$ intersect each other in two points which are not on any horizontal circle. See Fig.~\ref{pentagon}.  

Let us define the Alexander grading of a generator $\x \in \S(E)$ as follows: 

We can decompose a generator $\x$ of $C(E)$ as $\t\x_1\cup \t\x_2$ (non-uniquely) with the condition that the component of $\x$ on the circle $\gamma$ belongs to the $\t\x_1$. This is  possible since on each $\beta$ circle we have one component of $\x$ and on each $\alpha$ circle we have two components of $\x$ so we are in the situation of Lemma 3.1. of \cite{L}, and we can decompose $\x$ as the union of two generators, and we reorder them to satisfy the above condition.  Now consider $\t\x_1$ as an element of $C(H)$ and $\t\x_2$ as an element of $C(G)$, and define the Alexander grading of $\x$ as:
$$A(\x)=\frac{A_{C(H)}(\t\x_1)+A_{C(G)}(\t\x_2)}{2}$$

Here $A_{C(H)}(\t\x_1)$ means that we compute the Alexander grading with the markings of $H$, and similarly for $A_{C(G)}(\t\x_2)$ we use the markings of $G$. It is not hard to see that this is well-defined and is independent of the chosen decomposition.

We define a chain map (similar to \cite{MOST}) $\Phi_{\t \beta \gamma} : C(\t G)\longrightarrow C(E)$ by counting pentagons with the sign that we associate to each pentagon. Let $\x \in \S(\t G)$ and $\y \in \S(E)$, define $\pent $ to be the space of pentagons where each pentagon is an embedded disk in $\t T$ whose boundary consists of five arcs, each contained in a horizontal or a vertical circle, arranged as follows.
We start at the $\t \beta$ component of $\x$ and traverse the boundary with the orientation that comes from the pentagon, go through the arc of a horizontal circle $\A_i$, meet its corresponding $\y$ component, continue through an arc of a vertical circle $\B_k$, meet a component of $\x$, proceed to another horizontal circle $\A_j$, meet the component of $\y$ which is on the distinguished circle $\gamma$, continue along an arc in $\gamma$, meet an intersection point of $\t \beta$ with $\gamma$, which we call $a$. Finally traverse an arc in $\t \beta$ until we come back at the initial component of $\x$. All the angles here are required to be acute and all the pentagons are empty from $O$ or $X$ decorations and $\x$ components.
For later use, we represent a pentagon as above by the sequence $\A_i \rightarrow \B_k \rightarrow \A_j \rightarrow \gamma \rightarrow \t\B$ to keep track of the boundary arcs.
By using such pentagons, we connect $\x$ to $\y$, where the two generators differ from each other in exactly two components. See Fig.~\ref{pentagon}.

Fix an order for $\bm{\t\A}$ and $\bm{\t\B}$, also fix an orientation for each $\A$- and $\B$-circles. Let $\SS$ be a sign assignment of power $|\bm{\t\A}|$. Note that $\gamma$ has the same index in the ordering as $\B$. We orient $\gamma$ in the same way as $\B$. In this way the generators and flows in each diagram can be represented naturally by formal generators and formal flows.

In order to define the chain map, we need to assign signs to the pentagons. There are two types of pentagons: either the pentagon is on the left of the curve $\t \B$ which we call it a left pentagon, or on the right side of the curve $\t \B$ and we call it a right pentagon.

We associate to a pentagon $p\in \pent$ a formal rectangle $r(p)$. Let $p$ be represented by the sequence $\A_i \rightarrow \B_k \rightarrow \A_j \rightarrow \gamma \rightarrow \t\B$. If $p$ is a left (resp. right) pentagon then its associated formal rectangle is represented by the sequence $\A_i \rightarrow \B_k \rightarrow \A_j \rightarrow \t\B $ (resp. $\A_j \rightarrow \gamma \rightarrow \A_i \rightarrow \B_k $).

Define the sign for a pentagon as follows:

\begin{center}
$\varepsilon(p)=$
$\left \{ \begin{array}{ll}
\SS(r(p)) & \text{if $p$ is a left pentagon}   \\  
-\SS(r(p)) & \text{if $p$ is a right pentagon}
\end{array} \right. $ 
\end{center}

Given $\x \in \S(\t G)$, define

\begin{equation}
\Phi_{\t \beta \gamma}(\x) = \sum_{\y \in \S(E)} \sum_{p \in \pent} {\varepsilon(p)\y} \in C(E)
\end{equation}

Similarly we define $\Phi_{\gamma \t \B} : C(E)\longrightarrow C(\t G)$. Given $\x \in \S(E)$ and $\y \in \S(\t G)$, we define the space of pentagons $Pent_{\t \B \gamma}(\x,\y)$, where each pentagon is an embedded disk in $\t T$ whose boundary consists of five arcs, each contained in a horizontal or a vertical circle, arranged as follows.
With the boundary orientation, We start at the $\gamma$ component of $\x$, traverse the arc of a horizontal circle $\A_i$, reach its corresponding $\y$ component, traverse an arc of a vertical circle $\B_k$, meet a component of $\x$, continue through another horizontal circle $\A_j$, meet the $\t \B$ component of $\y$, continue along an arc in $\t \B$, meet an intersection point of $\t \beta$ with $\gamma$ (we call it $b$). At last, go through an arc in $\gamma$ to come back at the initial component of $\x$. All the angles here are required to be acute and all the pentagons are empty the components of $\x$ and $O$ or $X$ decorations.
We represent a pentagon as above by the sequence $\A_i \rightarrow \B_k \rightarrow \A_j \rightarrow \t\B \rightarrow \gamma$.
We denote by $r(p)$ the formal rectangle associated with $p$, which is represented by the sequence $\A_i \rightarrow \B_k \rightarrow \A_j \rightarrow \t\B$ if $p$ is a left pentagon, and with the sequence $\A_j \rightarrow \t\B \rightarrow \A_i \rightarrow \B_k$ if $p$ is a right pentagon.  We define the sign for the pentagon $p$ exactly the same as the previous case.


\begin{lemma}
The map $\Phi_{\t \beta \gamma} : C(\t G)\longrightarrow C(E)$ preserves the Alexander grading, and is an anti-chain map.
\end{lemma}

\begin{proof}
The fact that the Alexander grading will not change is straightforward.

The Argument to show that $\Phi_{\t \beta \gamma}$ is an anti-chain map is similar to Lemma 3.1 and Lemma 4.23 in \cite{MOST}. We consider different compositions of a pentagon and a rectangle, whether they are disjoint or have overlapping interiors or have an edge in common. In most cases the composite region has two different decompositions and the consistency of the signs follows from the property S-3 in Definition~\ref{def sign}. However there is one special case of composite regions that has a unique decomposition. They can be paired and each pair has the same initial and terminal points, one of the decompositions represents a term in $\partial \circ \Phi_{\t \beta \gamma}$ and the other one represents a term in $\Phi_{\t \beta \gamma} \circ \partial$; See Fig.~\ref{pentagon2s}.
Note that the formal rectangles associated with $r$ and $p$ form a $\B$-degeneration. Hence from the property S-2, we know that if $p$ is a left pentagon $\varepsilon(p) \cdot \SS(F(r))=-1$.
The formal rectangles associated with $r'$ and $p'$ make a $\B$-degeneration. Note that there is a minus in the definition of the sign of right pentagons so $\varepsilon(p') \cdot \SS(F(r'))=1$. Hence the two composite regions have opposite signs and cancel each other out.

\begin{figure}[h]
\centerline{\includegraphics[scale=0.7]{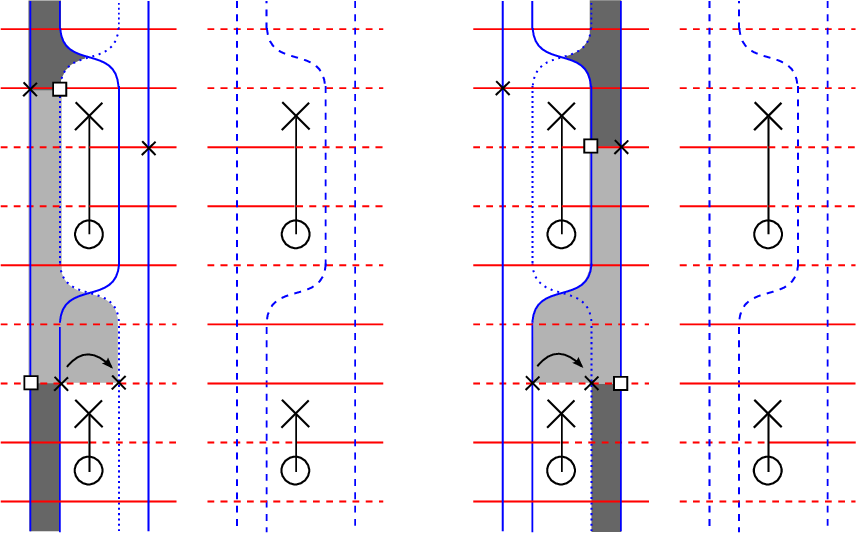}}
\caption {Here we showed a special case where two different composite regions have the same contribution. In the left diagram we first used the pentagon and then the rectangle. In the right diagram we first used a rectangle and then a pentagon. These two terms cancel out each other because the products of the signs are opposite.}
\label{pentagon2s}
\end{figure}

\end{proof}

\begin{figure}[h]
\centerline{\includegraphics[scale=0.7]{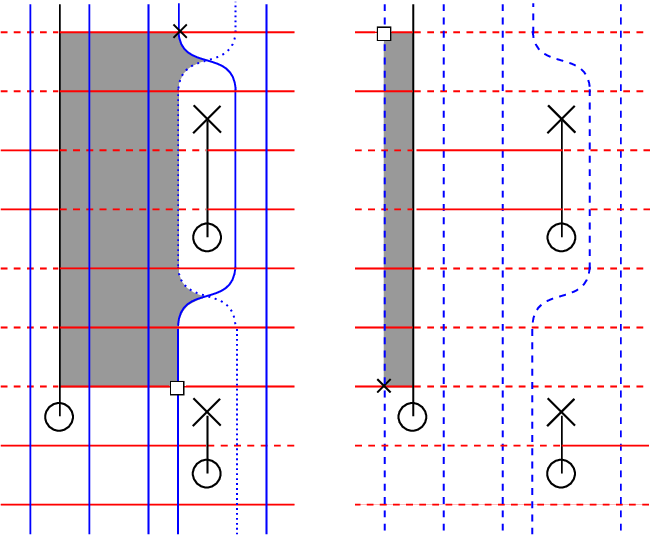}}
\caption {Here we showed a hexagon in $\hex$, where components of $\x$ are indicated by crosses and those of $\y$ are indicated by hollow squares.}
\label{hexagon}
\end{figure}

Following the ideas of \cite{MOST}, in order to define chain homotopy operators, we count hexagons. Given $\x, \y \in \S(\t G)$, we let $\hex$ denote the set of embedded hexagons in $\t T$. The boundary of a hexagon consists of six arcs, each contained in a horizontal or a vertical circle. More specifically, under the orientation induced on the boundary of $h$, we start at the $\t \beta$-component of $\x$, traverse the arc of a horizontal circle $\A_i$, meet its corresponding component of $\y$, continue through an arc of a vertical circle $\B_k$, meet its corresponding component of $\x$, proceed to another horizontal circle $\A_j$, meet its component of $\y$, which is contained in the distinguished circle $\t \beta$, continue along $\t \beta$ until the intersection point $b$ of $\t \beta$ and $\gamma$, proceed on $\gamma$ to the intersection point $a$ of $\t \beta$ with $\gamma$, and finally, continue on $\t \beta$ to the $\t \beta$-component of $\x$, which was also our initial point. All the angles of our hexagon are required to be acute and all the hexagons are empty from the $O$ or the $X$ decorations and the $\x$ components.
We represent a hexagon as above by the sequence $\A_i \rightarrow \B_k \rightarrow \A_j \rightarrow \t\B \rightarrow \gamma \rightarrow \t\B$.
See Fig.~\ref{hexagon}.

We now want to define signs for hexagons such that the maps defined by counting hexagons with respect to these signs are chain homotopy maps.
Let $h \in Hex_{\t \B \gamma \t \B} (\x , \y)$ be a hexagon that is represented by the sequence $\A_i \rightarrow \B_k \rightarrow \A_j \rightarrow \t\B \rightarrow \gamma \rightarrow \t\B$. We call $h$ a left (resp. right) hexagon if it is in the left (resp. right) side of $\t\B$. We associate to a left (resp. right) hexagon, the formal rectangle $r(h)$ that is represented by the sequence $\A_i \rightarrow \B_k \rightarrow \A_j \rightarrow \t\B$ (resp. $\A_j \rightarrow \t\B \rightarrow \A_i \rightarrow \B_k$).
We define the sign of $h$ to be $\varepsilon(h)=\SS(r(h))$. 

We define the homotopy operator  
$H_{\t \beta \gamma \t \beta}: C(\t G)\longrightarrow C(\t G)$
as follows:
\begin{equation}
H_{\t \beta \gamma \t \beta}(\x) = \sum_{\y \in \S(\t G)} \sum_{h \in \hex} \varepsilon(h) \y .
\end{equation}
Similarly we define $H_{\gamma \t \beta \gamma}: C(E)\longrightarrow C(E)$ over $\Z$.

\begin{proposition}
The map ${\Phi}_{\t \beta \gamma} : C(\t G) \longrightarrow C(E)$ is a chain homotopy equivalence with respect to sign assignments. In other words

\begin{center}
$\mathbb{I} + {\Phi}_{\gamma \t \beta} \circ {\Phi}_{\t \beta \gamma} + \partial \circ H_{\t \beta \gamma \t \beta} + H_{\t \beta \gamma \t \beta} \circ \partial = 0$

$\mathbb{I} + {\Phi}_{\t \beta \gamma} \circ {\Phi}_{\gamma \t \beta} + \partial \circ H_{\gamma \t \beta \gamma} +
H_{\gamma \t \beta \gamma} \circ \partial = 0.$
\end{center}

\end{proposition}

\begin{proof}

The proof is similar to the proof of \cite[Proposition 3.2 and 4.24]{MOST}.
Note that when we have a composite region, we consider the formal rectangle associated with the rectangle or pentagon or hexagon at hand. Hence the desired result came from the property S-3. 
\end{proof}

A similar argument about the Heegaard diagrams $\t H$ and $E$ proves the desired result.


\subsection{Stabilization} \label{stab}

Now that we proved commutation invariance, we want to show the stabilization invariance.
Let $\t G=(\t T, \bm{\t \alpha}, \bm{\t \beta}, \O, \X)$ be a Heegaard diagram and denote a stabilization of $\t G$ by $\t H=(\t U, \bm{\t \alpha '}, \bm{\t \beta '}, \O \cup O_1, \X \cup X_1)$ where $\bm{\t \A '} = \bm{\t \A} \cup \t\A^0_1 \cup \t\A^1_1$ and $\bm{\t \B '} = \bm{\t \B} \cup \t\B^0_1 \cup \t\B^1_1$.

Fix an ordering for the $\A$- and $\B$-circles in $\t H$ and orient them.
Let $\SS$ be a sign assignment of power $2n+2$. In this way, we can associate to each generator and flow of the diagram $\t H$ a formal generator and a formal flow of power $2n+2$.
We need to fix a sign assignment for the Heegaard diagram $\t G$. For each $\x\in S(\t G)$, by adding to $\x$ the two components $w_0$ and $w_1$ (with the intersection signs that come from the orientations of the corresponding arcs in $\t H$), we can also associate a formal generator $\x_f$ of power $2n+2$ to $\x$. Also to each flow $\phi:\x\rightarrow\y$ in $\t G$ we can associate a formal flow $F(\phi):\x_f\rightarrow\y_f$ of power $2n+2$.

Here we consider all complexes with coefficients in $\Z$ and the $\partial$ operator is defined with respect to the sign assignment in the end of Section~\ref{assign}. 
Note that the sign assignment of power $2n$ that we obtained from a sign assignment of power $2n+2$ for the complex $C(\t G)$ is gauge equivalent to any sign assignment of power $2n$. So it is well-defined to work with the above induced sign assignment.

Let $B = C(\t G)$ and $C=C(\t H)$ and $C'$
be the chain complex $B[1]\oplus B$ ($B[1]$ is the chain complex obtained from $B$ by shifting the Alexander grading by $1$), endowed with the differential
 $\partial' : C'\longrightarrow C'$ given by
\[
\partial'(a,b)=(\partial a,-\partial b),
\]
where $\partial$ denotes the differential within $B$. Note that $C'$ is the mapping cone of the zero map between $B$ and itself. Let $\L$ and $\RR\cong B$ be the subgroups of $C'$ of elements of the form $(c,0)$ and $(0,c)$ for $c\in B$, respectively. The module $\RR$ inherits the Alexander grading from its identification with $B$ and $\L$ is given the Alexander grading which is one less than the one it inherits from its identification with $B$. This shows that $H(C')=H(B)\otimes V$ where $V\cong \F \oplus \F$ with generators in gradings $0$ and $-1$, where $\F$ is the ring of coefficient; in this case $\Z$.

Consider the stabilized grid diagram $H$ for the knot $K$. We denote by $X_2$, the element of $\X$ which is placed in the row that $O_1$ is placed.
As we mentioned in the definition of stabilization move, we only need to consider stabilizations such that $X_2$ is placed in the square just to the left or just to the right of $O_1$ (in the grid diagram $H$).

We denote by $w_0$ (resp. $w_1$), the intersection of $\t \B_1^0$ (resp. $\t \B_1^1$) with the lifts of $\A_1$. Let $(I,I)\subset \S(\t H)$ be the set of those generators which have both of their $\t\alpha^i_1$ components for $i=0,1$ on 
the lifts of $\beta_1$ (i.e. the generators which have $w_0$ and $w_1$ as its components.) There is a natural one-to-one correspondence between elements of $\S(\t G)$ and elements of $(I,I)$ with the following property: For $\x \in \S(\t G)$, let $\psi(\x) \in \S(\t H)$ be the associated generator in $(I,I)$. We have 
\begin{equation}
A_{C(\t G)}(\x)=A_{C(\t H)}(\psi (\x))+1=A_{C'}(0,\x)=A_{C'}(\x,0)+1
\end{equation}

In order to define a map from $C(\t H)$ to either $\L$ or $\RR$ we need to define some combinatorial objects.

Let $\gamma$ be an arc that connects $O_1$ to $X_1$ in the Heegaard diagram $H$ for $K$. Denote by $\t \gamma$ the preimage of $\gamma$ in the the Heegaard diagram $\t H$ for $\t K$. Note that $\t \gamma$ is a circle.

\begin{definition}
Let $\x,\y\in \S(\t H)$, a \emph{pseudo-domain} from $\x$ to $\y$ is a two-chain in the Heegaard diagram $\t H$ whose boundary consists of a path from $\x$ to $\y$ and possibly a number of copies of $\gamma$. We denote by $\sigma(\x,\y)$ the set of pseudo-domains from $\x$ to $\y$. Note that $\pi(\x,\y) \subset \sigma(\x,\y)$.
\end{definition}

\begin{remark} In this definition we allow the two-chains to have a number of copies of $\gamma$ in their  boundary, but for most parts of the paper we only work with regions with at most one such copy. The reason that we allow multiples of $\gamma$ is to be able to extend the $*$ operator to $\sigma(\x,\y)$.
\end{remark}

We introduced the $*$ operation earlier. There is a natural extension of $*$ to pseudo-domains.
Given $\mathbf{a} , \mathbf{b} ,\mathbf{c} \in S(\t H)$, if $p \in \sigma(\mathbf{a}, \mathbf{b})$ and $p' \in \sigma (\mathbf{b} , \mathbf{c})$, we can add them as two-chains and by the definition of pseudo-domains it will be an element of $\sigma(\mathbf{a},\mathbf{c})$. In this way we get
$$*:\sigma(\mathbf{a}, \mathbf{b})\times \sigma (\mathbf{b} , \mathbf{c})\longrightarrow \sigma(\mathbf{a},\mathbf{c}) $$

\begin{definition}\label{prect} Let $\x, \y \in S(\t H)$ be two generators that differ exactly along $\t\B_1^0$ and $\t\B_2^0$. A \emph{punctured rectangle} $\a$ that connects $\x$ to $\y$ is a topologically embedded punctured disk in $\t U$ (the Heegaard surface of $\t H$). The boundary of $\a$ consists of the circle $\gamma$ and a path from $\x$ to $\y$ consisting of four arcs such that $\partial \a \cap \bm{\t \beta '} \subset \t\B_1^0 \cup \t\B_2^0$ or $\partial \a \cap \bm{\t \beta '} \subset \t\B_1^1 \cup \t\B_2^1$. We denote by $A(\x , \y)$ the set of punctured rectangles from $\x$ to $\y$.
\end{definition}

Although $\a \in A(\x , \y)$ is not a domain, for our combinatorial purposes we need to connect the generators as above. Note that $A(\x , \y)$ has at most 1 element.

Given $\x , \y \in S(\t H)$ that differ exactly along $\t\B_1^0$ and $\t\B_2^0$, let $\mathfrak{a} \in A(\x, \y)$ be a punctured rectangle. There is a unique empty rectangle in $Rect(\y, \x)$ that we denote by $r_{\mathfrak{a}}$ and we call it the complementary rectangle of $\a$. Note that the support of the union of ${\mathfrak{a}}$ and $r_{\mathfrak{a}}$ is topologically a sphere with three punctures. Define $\mu(\a) := -\SS(F(r_{\a}))$, where $F(r_{\a})$ is the formal rectangle corresponding to $r_{\a}$. 

There is another way to define the sign of ${\mathfrak{a}}$. We associate to ${\mathfrak{a}}$, the formal rectangle $F({\mathfrak{a}}): \x_f \rightarrow \y_f$ in $\FF_{2n+2}$ that has the same boundary arcs as ${\mathfrak{a}}$, where $\x_f$ is the formal generator associated with $\x$. Define $\mu(\a) := \SS(F({\a}))$. This definition gives the same sign for ${\mathfrak{a}}$ as the previous definition, because of the property S-2 of the definition of sign assignments.

\begin{definition} \label{lshape} For $\x \in S(\t H)$, $w_i \in \y \in S(\t H)$, we define an \emph{L-shape} associated with $w_i$ ($i= 0, 1$) that connects $\x$ to $\y$ to be an embedded punctured disk that is empty from the components of $\x$ and basepoints except for $O_1$ and $X_1$.
Its boundary consists of the circle $\gamma$ and a path from $\x$ to $\y$ that consists of six arcs, arranged as follows. We start from $w_i$ which is a component of $\y$ and traverse the boundary with the orientation that comes from the L-shape (counterclockwise).
We go through $\t \B _1^i$ we reach at a component of $\x$. Then we traverse an $\A$-curve to reach the component of $\y$ on $\t\B_2^i$, then going through $\t\B_2^i$ we meet a component of $\x$. Then we traverse an $\A$-curve, namely $\A_j$, to meet a component of $\y$, we go through a $\B$-curve, namely $\B_k$, and reach a component of $\x$. Finally we go along a lift of $\A_1$ to reach $w_i$. See Fig.~\ref{L}. Alternately we can define an L-shape $l$ to be the composite region made
up of a punctured rectangle $a \in A(\x , \z)$
and $r\in Rect(\z, \y)$, where $\z \in S(\t H)$ differs with $\x$
exactly along $\t \B_1^i$ and $\t \B_2^i$, and the upper
(resp. lower) edge of $r$ is an arc in $\A_j$ (resp. $\A_1$), and the
right (resp. left) edge of $r$ is an arc in $\t \B_1^i$ (resp. $\B_k$), i.e. $l = \a *r$.

We define the sign of an L-shape $l = \a *r$ (associated with $w_0$ or $w_1$) that connects $\x$ to $\y$ as follows:
$$\mu(l) := \mu(\a) \cdot \SS(F(r)),$$
\noindent where $F(r)$ is the formal rectangle associate with the empty rectangle $r$. 
\end{definition}

\begin{definition} \label{Ldomain} Let $\x \in S(\t H)$ and $\y \in (I,I) \subset S(\t H)$. A pseudo-domain $p\in \sigma (\x,\y)$ is called \emph{Type L} if one of the following is true:

\textbf{L(1)} $\x=\y$, and $p$ is the trivial domain. In this case we define the sign of $p$ as $\mu(p):= 1$.

\textbf{L(2)} $\x$ contains either $w_0$ or $w_1$ but not both, and $p$ is an L-shape. See Fig.~\ref{L}. In this case, the sign of $p$ is $\mu(p)$ as defined above.

\textbf{L(3)} $\x$ contains neither $w_0$ nor $w_1$, and $p$ is the sum of two L-shapes $l_0$ and $l_1$. Note that in this case the branched cuts inside the L-shape regions are glued together to form a domain whose boundary consists of twelve arcs. In this case, we define $\mu(p):= \mu(l_0)\cdot \mu(l_1)$.

\begin{figure}[h]
\centerline{\includegraphics[scale=0.8]{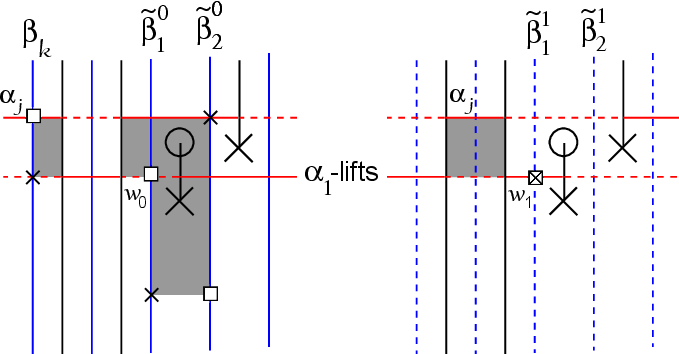}}
\caption {We have illustrated a pseudo-domain of Type L which consists of one L-shape and $w_1 \in \x$. Note that the branch cuts are shown in black.}
\label{L}
\end{figure}

\end{definition}

\begin{figure}[h]
\centerline{\includegraphics[scale=0.9]{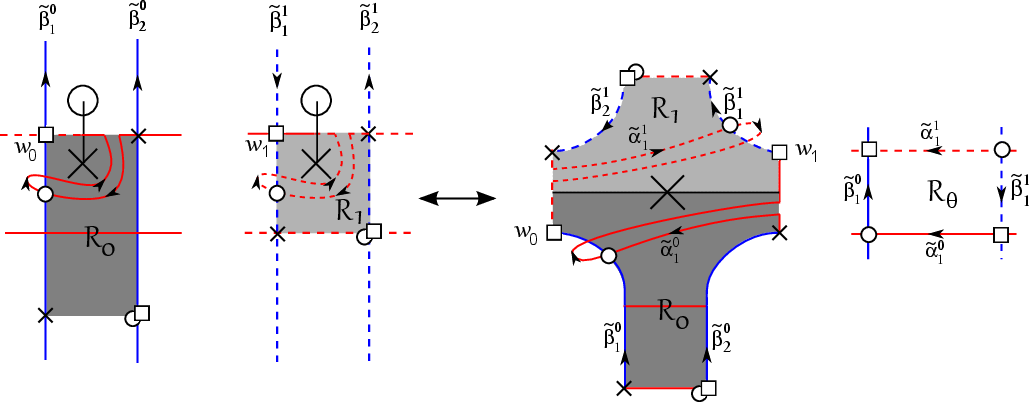}}
\caption {Here we have illustrated an octagon and the isotopies that we used in the definition of its sign. Also we showed the orientations of the edges of the octagon, and according to that we oriented the edges of the formal rectangle $R_{\theta}$. The generator $\x$ is shown by crosses. The generators $\w$ and $\y$ (and their associated formal generators) are shown by hollow circles and hollow squares, respectively.
}
\label{octagon}
\end{figure}

\begin{definition} \label{oct} For $\x \in S(\t H)$ and $\y \in (I,I) \subset S(\t H)$, an \emph{octagon} $\theta \in \pi(\x , \y)$ is topologically an embedded disk in $\t U$ (the Heegaard surface of $\t H$) whose boundary is a path form $\x$ to $\y$ that consists of eight arcs, arranged as follows. Starting from the component of $\x$ on $\t\B_1^0$, we traverse the boundary with the orientation that comes from the octagon. We go along an $\A$-curve to meet a $\y$ component. We traverse a $\B$-curve to reach a component of $\x$, then we go through one of the lifts of $\A_1$, passing the branch cut connecting $O_1$ and $X_1$, we meet a component of $\y$ that is on the $\t\B_1^1$ and is the same as $w_1$. We go through $\t\B_1^1$ and get to a component of $\x$, after traversing an $\A$-curve we reach a $\y$ component, then we go along a $\B$-curve to meet one of $\x$ components. Going through the other lifts of $\A_1$, we pass the branched cut connecting $O_1$ and $X_1$ and we reach a component of $\y$ which is on $\t\B_1^0$ and is equal to $w_0$. Finally we go along $\t\B_1^0$ and come back to the $\x$ component that we started from. The interior of $\theta$ is empty from the components of $\x$ and basepoints other that $X_1$. See Figs.~\ref{R} and \ref{octagon}. 

\begin{figure}[h]
\centerline{\includegraphics[scale=0.8]{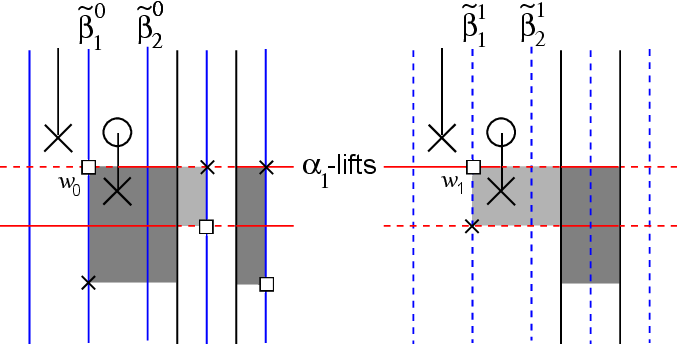}}
\caption {An octagon that is obtained by gluing together the two parts with different shadings. This octagon connects the generator $\x$ which is shown by crosses to the generator $\y$ which is shown by hollow squares.}
\label{R}
\end{figure}

\end{definition}

Given $\x \in S(\t H)$ and $\y \in (I,I) \subset S(\t H)$, let $\theta \in \pi(\x , \y)$ be an octagon.
Consider an arc $\gamma_i$ in the $i^{th}$ grid ($i=0,1$). The arc starts from a point on the lift of $\A_1$ just to the right of the branched cut connecting $O_1$ and $X_1$.
The arc $\gamma_i$ intersects $\t\B_1^i$ and goes to the square (or octagon) just to the left of $X_1$ in the $i^{th}$ grid.
We isotope each lift of $\A_1$ by doing a finger move along the appropriate arc, either $\gamma_0$ or $\gamma_1$. See Fig.~\ref{octagon}.

Without loss of generality, we assume that $\gamma_0$ (resp. $\gamma_1$) starts from a point on $\t\A_1^0$ (resp. $\t\A_1^1$).
Let $R_0 \in Rect(\x , \z)$ be the new rectangle that is created after the isotopy of $\t\A_1^0$ along $\gamma_0$, and let $R_1 \in Rect(\z , \w)$ be the other rectangle which is created after the isotopy along $\gamma_1$, where $\z$ and $\w$ are generators in $S(\t H)$. See Fig.~\ref{octagon}. With the orientations that we fixed for $\bm \A$ and $\bm \B$ in the beginning of this section, we can associate to each generator $\x \in \S(\t H)$ and each empty rectangle $r$, a formal generator $\x_f$ and a formal rectangle $F(r)$.

Denote by $R_{\theta}$ the formal rectangle from $\w_f$ to $\y_f$, such that the sequence $\t\A_1^1 \rightarrow \t\B_1^0 \rightarrow \t\A_1^0 \rightarrow \t\B_1^1$ represents its edges (with the notation of Definition~\ref{formalrec}), and we orient its edges so that the intersection numbers of the lower-left and upper-right corners coincide with the intersection numbers of the components of $\w_f$, and the intersection numbers of the lower-right and upper-left corners coincide with the intersection numbers of the components of $\y_f$. See Fig.~\ref{octagon}. 

We define the sign for the octagon $\theta$ as follows:
\[ \mu(\theta):=\SS(F(R_0))\cdot \SS(F(R_1))\cdot \SS(R_{\theta}), \]

\begin{definition} \label{Rdomain} Let $\x \in S(\t H)$ and $\y \in (I,I) \subset S(\t H)$. A pseudo-domain $p\in \sigma(\x,\y)$ is called \emph{Type R} if one of the following is true:

\textbf{R(1)} There exists an octagon $\theta \in \pi(\x, \y)$, and $p=\theta$. In this case, the sign $\mu(p)$ is defined as above.

\textbf{R(2)} There are intermediate generators $\mathbf{u} , \mathbf{v} \in S(\t H)$ such
that there exist $\a\in A(\x , \mathbf{ u})$, $r\in Rect(\mathbf{ u} , \mathbf{v})$ and an octagon $\theta \in \pi (\mathbf{ v} , \y)$ with following properties for $i=0$ or $i=1$:
\begin{itemize}
\item $\mathbf{u}$ and $\x$ differ exactly along $\t \B_1^i$ and $\t \B_2^i$.
\item the right edge of $r$ is contained in $\t \B_1^i$.
\item the upper (resp. lower) edge of $r$ is contained in the same $\A$-circle as the upper edge of $\a$ (resp. as one of $\theta$ edges).
\item the lower-right corner of $r$ which is the component of $\mathbf{v}$ along $\t \B_1^i$, should be contained in the segment of the left edge of $\a$ between $w_i$ and the lower-left corner of $\a$.
\end{itemize}
In this case $p=\a *r*\theta$ and we define its sign to be $\mu(p) : = \mu(\a) \cdot \SS(F(r))\cdot \mu(\theta)$. See Fig.~\ref{R(2)}.

\textbf{R(3)}\hspace{5pt} There are intermediate generators $\mathbf{ u} , \mathbf{ u'} ,\mathbf{v} , \mathbf{ v'} \in S(\t H)$ such that there exist $\a \in A(\x , \mathbf{ u})$, $r\in Rect(\mathbf{u} , \mathbf{v})$, $\a '\in A(\mathbf{v} , \mathbf{u'})$, $r'\in Rect(\mathbf{u'} , \mathbf{v'})$ and an octagon $\theta \in \pi(\mathbf{v'} , \y)$ with following properties:
\begin{itemize}
\item $\x$ and $\mathbf{u}$ differ exactly along $\t \B_1^0$ and $\t \B_2^0$.
\item $\mathbf{v}$ and $\mathbf{u'}$ differ exactly along $\t \B_1^1$ and $\t \B_2^1$.
\item the right edge of $r$ (resp. $r'$) is contained in $\t \B_1^0$ (resp. $\t \B_1^1$).
\item the upper (resp. lower) edge of $r$ is contained in the same $\A$-circle as the upper edge of $\a$ (resp. one of $\theta$ edges).
\item the upper (resp. lower) edge of $r'$ is contained in the same $\A$-circle as the upper edge of $\a '$ (resp. one of $\theta$ edges).
\item the lower-right corner of $r$ which is the component of $\mathbf{v}$ along $\t \B_1^0$, should be contained in the segment of the left edge of $\a$ between $w_0$ and the lower-left corner of $\a$.
\item the lower-right corner of $r'$ which is the component of $\mathbf{v'}$ along $\t \B_1^1$, should be contained in the segment of the left edge of $\a '$ between $w_1$ and the lower-left corner of $\a '$.
\end{itemize}
In this case $p=\a *r* \a '*r'*\theta$ and we define its sign to be $\mu(p) : = \mu(\a) \cdot \SS(r) \cdot \mu(\a ') \cdot \SS(F(r')) \cdot \mu(\theta)$.

\begin{figure}[h]
\centerline{\includegraphics[scale=0.8]{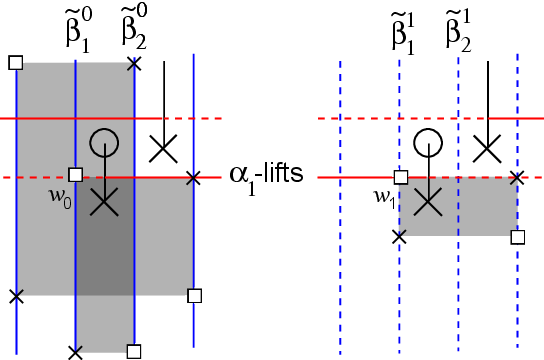}}
\caption {We have illustrated a pseudo-domain of the form R(2) that connects the generator shown by crosses to the generator shown by hollow squares.}
\label{R(2)}
\end{figure}

\end{definition}

We now define the maps 
\[F^L:C(\t H)\rightarrow\mathcal{L}\]
\[F^R:C(\t H)\rightarrow\mathcal{R},\]
where $F^L$ (resp. $F^R$) counts pseudo-domains of Type $L$ (resp. $R$), more precisely define
\begin{eqnarray*}
F^L(\x)= \sum_{\y \in \S(\t H)}\sum_{p\in \sigma^{L}(\x,\y)} \mu(p) \y \\
F^R(\x)= \sum_{\y \in \S(\t H)}\sum_{p\in \sigma^{R}(\x,\y)} \mu(p) \y ,
\end{eqnarray*}
\noindent where $\sigma^{L}(\x,\y)$ (resp. $\sigma^{R}(\x,\y)$) denotes the set of pseudo-domains of Type $L$ (resp. $R$). We set $\sigma^{F}(\x,\y)=\sigma^{L}(\x,\y)\cup \sigma^{R}(\x,\y)$ and define 

\[F=\begin{pmatrix} F^L\\ F^R \end{pmatrix}:C(\t H)\rightarrow C'\]

Now we want to show that $F$ is a chain map. Note that in the boundary map that appears in $\partial\circ F$ we count empty rectangles in $\t G$. Recall that there is a one-to-one correspondence between $\S(\t G)$ and the set $(I,I)\subset \S(\t H)$ that is given by $\psi$.

\begin{figure}[h]
\centerline{\includegraphics[scale=0.8]{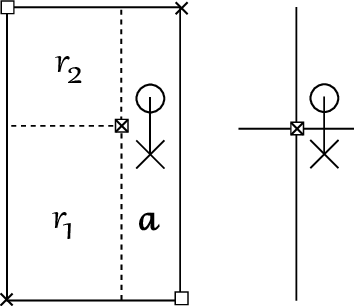}}
\caption {We have illustrated a pseudo-rectangle of Type 2, that can be counted in boundary map $C(\t G)$ which is drawn in $\t H$.}
\label{Type2rect}
\end{figure}

\begin{definition} Given $\x , \y \in \S(\t G)$ and an empty rectangle $r\in Rect(\x,\y)$, there is a pseudo-domain $r'$ between the two associated generators $\psi(\x),\psi(\y)\in(I,I)$ (not necessarily empty). By abuse of notation we denote $r'$ by $\psi(r)$. If $\psi(r)$ is a topological rectangle in $\t H$ we call it a Type 1 pseudo-rectangle. On the other hand there are empty rectangles $r$ in $\t G$ that $\psi(r)$ is an annuli that contains either $w_0$ or $w_1$ in its interior. Its boundary (aside from the circle $\gamma$) is a path from $\psi(\x)$ to $\psi(\y)$ which is made of four arcs contained in $\A$- and $\B$-circles other than the circles $\t \A_1^i$ or $\t \B_1^i$ for $i=0, 1$. We call these pseudo-domains Type 2 pseudo-rectangles (See Fig.~\ref{Type2rect}).
\end{definition}

The boundary map of $C(\t G)$ can be computed by counting pseudo-rectangles in $\t H$. Note that a Type 2 pseudo-rectangle is a pseudo-domain in $\t H$ of the form $r_1 * \a * r_2$, where $r_1 \in Rect (\psi(\x) , \mathbf{u})$, $\a \in A(\mathbf{u} , \mathbf{v})$ and $r_2 \in Rect (\mathbf{v} , \psi(\y))$. Here $\mathbf{u} , \mathbf{v} \in \S(\t H)$ are intermediate generators. Denote by $F(r)$ the formal rectangle that has the same boundary arcs as $r$ (which is the same as the boundary arcs of $\psi(r)$). We define the sign of $\psi(r)$ to be the sign of $F(r)$, i.e. $\SS(r)=\SS(F(r))$.


\begin{remark} \label{rem} Let $\x_f$ and $\y_f$ be formal generators of power $m$ and $\phi:\x_f \rightarrow \y_f$ be a formal rectangle that is embedded in a grid diagram with the support $ABCD$, where $A$, $B$, $C$ and $D$ meet in a component of the initial generator in the interior. See Fig.~\ref{pointedrect}. The proof of Proposition~4.4 in~\cite{Sign} shows that the sign of formal rectangles satisfies the following identities:
\begin{eqnarray*}
\SS(ABCD)
&=&\SS(A)\cdot \SS(BC)\cdot \SS(D)\\
&=&\SS(C)\cdot \SS(AD)\cdot \SS(B),
\end{eqnarray*}
\noindent where the alphabets show the support of the formal rectangles with appropriate initial generators.
For example, in the first equality $A$ denotes the formal rectangle $r_A :\x_f \rightarrow \z_f$, $BC$ is the support of the embedded rectangle associated with $r_{BC}: \z_f \rightarrow \w_f$ and $r_{D} :\w_f  \rightarrow \y_f$ is denoted by $D$, where $\z_f$ and $\w_f$ are the formal generators of power $m$.
In this case the formal rectangle $\phi$ can alternatively described as the composition of the three formal rectangles $r_A$, $r_{BC}$ and $r_{D}$.

\begin{figure}[h]
\centerline{\includegraphics[scale=0.9]{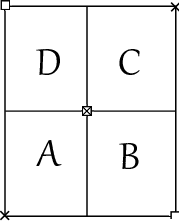}}
\caption {Here we have illustrated a formal rectangle that is embedded in a grid diagram as a rectangle with a component of the initial generator in its interior and with the support $ABCD$. The initial generator $\x_f$ is shown with crosses. The terminal generator $\y_f$ is shown with hollow squares.}
\label{pointedrect}
\end{figure}

\end{remark}

\begin{lemma} Let $r$ be a Type~2 pseudo-rectangle that is defined by $r_1 * \a * r_2$. Then
$$\SS(r)=\SS(r_1)\cdot \mu(\a)\cdot \SS(r_2).$$
\end{lemma}

\begin{proof} The proof follows from~\cite[Proposition~4.4]{Sign} as explained in Remark~\ref{rem}.
\end{proof}

We consider only those stabilizations, that the basepoint $X_2$ is placed just to the left or just to the right of $O_1$ in the grid diagram $H$. If $X_2$ is placed just to the left of $O_1$, by definition, we can not have any L-shapes, Type~2 pseudo-rectangles, or pseudo-domains of the form R(2) and R(3).
However, if $X_2$ is placed just to the right of $O_1$, since the pseudo-domains can not contain $X_2$, there are some restrictions on them.

Given $\x ,\y \in (I,I)\subset S(\t H)$, let $r\in Rect(\x, \y)$ be an empty rectangle in $\t G$. Note that the circles $\t\A_1^i$ and $\t\B_1^i$ (for $i=0,1$) are not in the sets of $\A$ and $\B$-circles for the Heegaard diagram $\t G$ (these circles are added in the stabilization move). Hence the arcs in the boundary of $r$ can not be contained in $\t\A_1^i$ and $\t\B_1^i$ for $i=0,1$. We first show that $F$ is a chain map over $\Z_2$ and later we prove it with $\Z$ coefficients. 

\begin{lemma} \label{FZ2}
The map $F:C(\t H)\rightarrow C'$ preserves the Alexander grading, and is a chain map over $\Z_2$.
\end{lemma}

The idea of the proof is similar to the proof of~\cite[Lemma 3.5]{MOST}, but we give a complete proof and discuss all the possible configurations and include all the necessary changes to their proof.

\begin{proof}

\begin{figure}[!h]
\centerline{\includegraphics[scale=0.7]{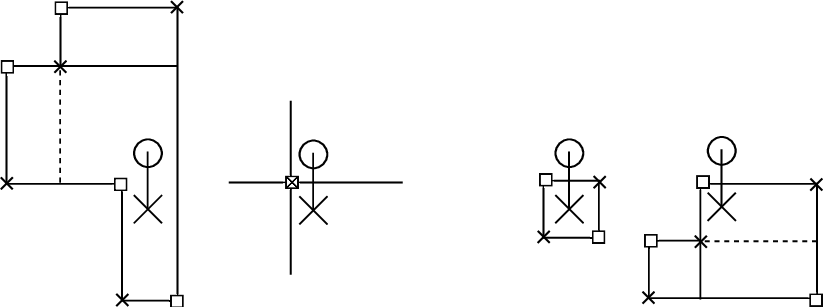}}
\caption {\textbf{Contributions from I(1)}
The first Figure shows a Type $L$ pseudo-domain and juxtaposing rectangle. The second Figure shows an octagon and a juxtaposing rectangle.}
\label{I1}
\end{figure}

\begin{figure}[h!]
\centerline{\includegraphics[scale=0.9]{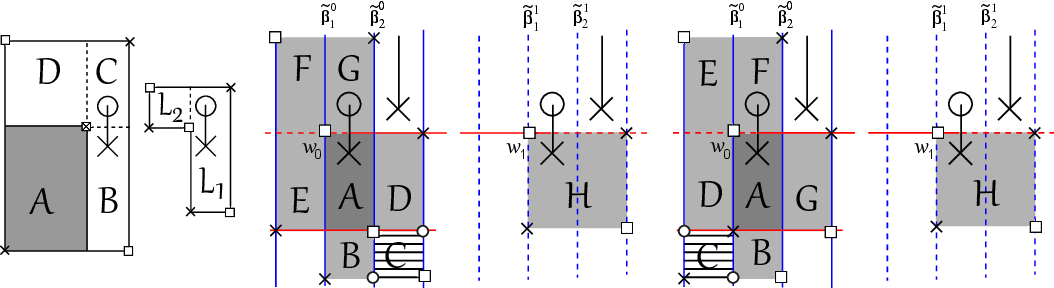}}
\caption {\textbf{Contributions from I(2) and II(0).} In the above pictures, the initial generator $\x$ is shown with crosses. The terminal generator $\y$ is shown with hollow squares. The hollow circles show the intermediate generators.
The left picture shows a term $p*r$ from II(0), where $p=L_1L_2$ is a Type L pseudo-domain of the form L(2) and $r=ABCD$ is a Type 2 pseudo-rectangle, which cancels out with a term $r'*p'$ from I(2), where $r'=A$ is a Type 1 pseudo-rectangle and $p'=BCDL_1L_2$ is a Type L pseudo-domain of the form L(3).
The next picture shows a contribution $p*r$ from II(0) that cancels out with a contribution $p'*r'$ from I(2), where $p, p'\in \sigma^R$, $p=ABCDH$ is an octagon and $p'=ABDEFGH$ is of the form R(2), $r'=C$ is a Type 1 pseudo-rectangle (shaded horizontally) and $r=AEFG$ is a Type 2 pseudo-rectangle.
The right picture shows a contribution $p*r$ from II(0) that cancels out with a contribution $r'*p'$ from I(2), where $p, p'\in \sigma^R$, $p=AGH$ is an octagon and $p'=ABDEFGH$ is of the form R(2), $r'=C$ is a Type 1 pseudo-rectangle (shaded horizontally) and $r=ABCDEF$ is a Type 2 pseudo-rectangle.}
\label{newII0}
\end{figure}

The fact that the Alexander grading is preserved, comes from the definition of Alexander grading and the way we have defined $\mathcal{L}$ and $\mathcal{R}$ with appropriate gradings.

We group together different possibilities for the terms in $\partial \circ F$ and $F\circ\partial$. Using basic planar geometry, we arrange the cases of the composition of an empty rectangle $r$ and a pseudo-domain $p\in \sigma^F$, according to the type of the empty pseudo-rectangle $r$, and the number of common corners of $r$ and $p$. Note that the image of the support of the composite of $r$ and $p$ in the grid diagram $H$, can not wrap horizontally around the Heegaard surface associated with $H$ (i.e. the torus). Because this composition can not contain basepoints (other than possibly $O_1$ and $X_1$).

If $r$ is of Type 1, we have the following cases:

I(0) A composition of a pseudo-rectangle $r$ of Type 1 and a pseudo-domain $p\in \sigma^F$, where they do not have any corners in common. This composition can be counted in either of $\partial \circ F$ or $F\circ\partial$. 

I(1) A composition of $r$ (a pseudo-rectangle of Type 1) and $p\in \sigma^F$ in either order, where they have one corner in common and $r$ does not contain $w_0$ (or $w_1$) in its boundary. This composition has a unique concave corner, cutting through this concave corner in either of the two possible manners results in a decomposition. See Fig.~\ref{I1}.

I($1'$) A composition of a pseudo-rectangle $r$ of Type 1 and $p\in \sigma^F$, where they share one corner and $w_0$ (or $w_1$) is in the boundary of $r$. In this case, at least one of the edges of $r$ is contained in the circles that are not in $\t G$, hence the composite region is of the form $r*p$ and is counted in $F\circ \partial$. Considering that $r$ is empty from basepoints, there are two possibilities here: First $w_0$ (or $w_1$) is in the right edge of $r$. Second, the lower-right corner of $r$ is $w_0$ (or $w_1$).

I(2) A composition in either order of a pseudo-rectangle $r$ of Type 1 and $p\in \sigma^F$, where they share two corners other than possibly $w_0$ (or $w_1$). See Fig.~\ref{newII0}.

I(3) A composition of a pseudo-rectangle $r$ of Type 1 and $p\in \sigma^F$, where they share three corners other than possibly $w_0$ (or $w_1$). In this case we can see that the composite region contains one of the two columns that contain $O_1$ and $X_1$. See Figs.~\ref{newI3L} and \ref{newI3R}.

\begin{figure}[h!]
\centerline{\includegraphics[scale=0.8]{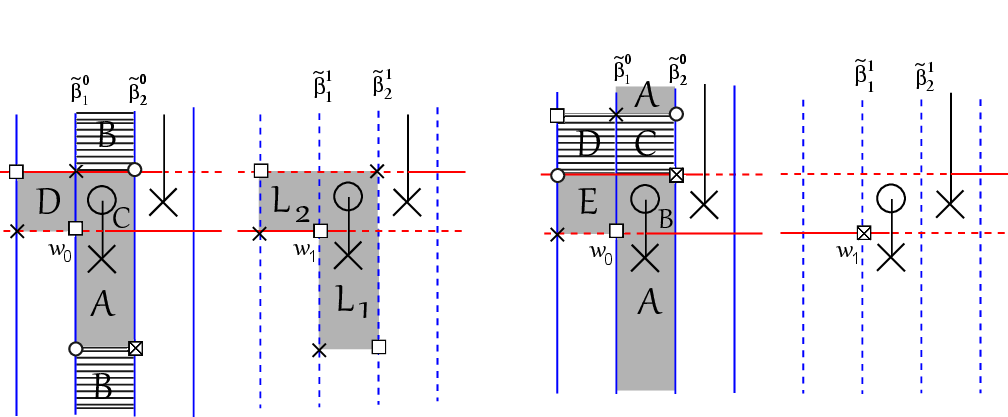}}
\caption {\textbf{Contributions from I(3) and I($1'$)}
The left picture shows a contribution $r*p$ from I(3) that cancels out with a contribution $r'*p'$ from I($1'$), where $p, p'\in \sigma^L$, $p=ACDL_1L_2$ is of the form L(3) and $p'=L_1L_2$ is of the form L(2), $r=B$ is shaded horizontally and $r'=D$ has $w_0$ as its lower-right corner.
The right picture shows a contribution $p*r$ from I(3) that cancels out with a contribution $r'*p'$ from I($1'$), where $p, p'\in \sigma^L$, $p=ABE$ is of the form L(2) and $p'$ is trivial, $r=CD$ is shaded horizontally and $r'=DE$ has $w_0$ as its lower-right corner.}
\label{newI3L}
\end{figure}

\begin{figure}[h!]
\centerline{\includegraphics[scale=0.8]{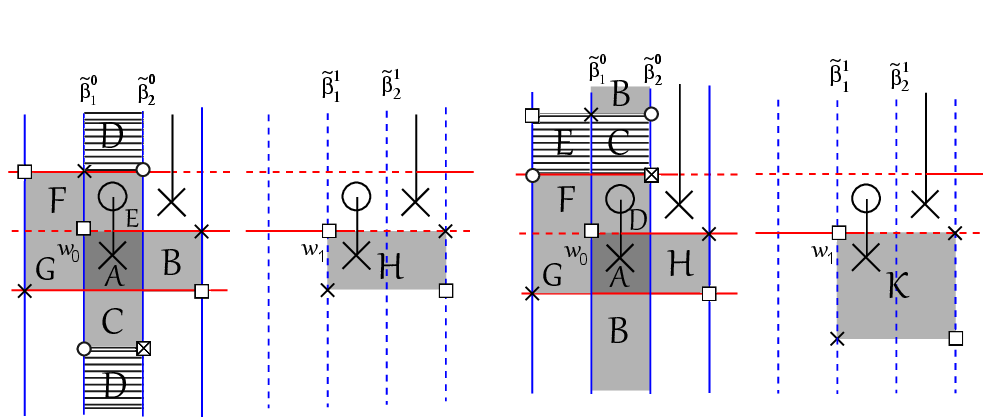}}
\caption {\textbf{Contributions from I(3) and I($1'$)}
The left picture shows a contribution $r*p$ from I(3) that cancels out with a contribution $r'*p'$ from I($1'$), where $p, p'\in \sigma^R$, $p=ABCEFGH$ is of the form R(2) and $p'=ABH$ is an octagon, $r=D$ is shaded horizontally and $r'=FG$ has $w_0$ in the interior of its right edge.
The right picture shows a contribution $p*r$ from I(3) that cancels out with a contribution $r'*p'$ from I($1'$), where $p, p'\in \sigma^R$, $p=ABDFGHK$ is of the form R(2) and $p'=AHK$ is an octagon, $r=CE$ is shaded horizontally and $r'=EFG$ has $w_0$ in the interior of its right edge.}
\label{newI3R}
\end{figure}

If $r$ is of Type 2, it can only appear in the terms of the form $\partial\circ F$ and we have the following cases:

II(0) The pseudo-rectangle $r$ of Type 2 has no corner in common with $p\in \sigma^F$, and the composite region is of the form $p*r$. See Fig.~\ref{newII0}.

II(1) A composite region of the form $p*r$, where $r$ is a pseudo-rectangle of Type 2 that shares one corner with $p\in \sigma^F$. We can see that in this case the composite region contains one of the two columns that contain $O_1$ and $X_1$. See Fig.~\ref{newII1}.

\begin{figure}[h!]
\centerline{\includegraphics[scale=0.8]{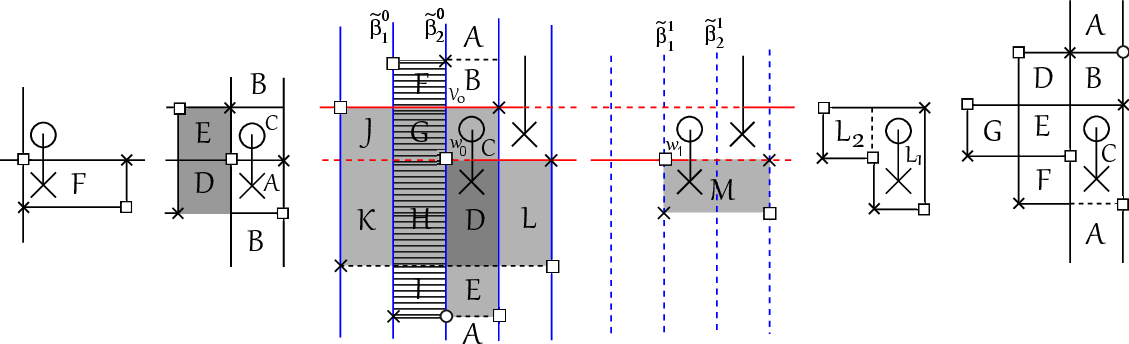}}
\caption {\textbf{Contributions from II(1) and I(1$'$)}
The left picture shows a term $p*r$ from II(1) that cancels out a term $r'*p'$ from I(1$'$), where $p, p' \in \sigma^R$, $p=ABF$ and $p'=AF$ are an octagons, $r=ACDE$ is a Type 2 pseudo-rectangle, and $r'=DE$ is shaded gray.
The next picture shows a term $p*r$ from II(1) that cancels out a term $r'*p'$ from I(1$'$), where $p, p' \in \sigma^R$, $p=ACDEGHJKLM$ is of the form R(2), $p'=CDEGHJKLM$ is of the form R(2) and is in gray, $r=BCDEFGHI$ is a Type 2 pseudo-rectangle, and $r'=FGHI$ is shaded horizontally.
The right picture shows a term $p*r$ from II(1) that cancels out a term $r'*p'$ from I(1$'$), where $p, p'\in \sigma^L$ are of the form L(3), $p=ACEGL_1L_2$ and $p'=CEGL_1L_2$. $r=BCDEF$ is a Type 2 pseudo-rectangle, and $r'=DEF$ is a Type 1 pseudo-rectangle that has $w_1$ in the interior of its right edge.}
\label{newII1}
\end{figure}

Contributions from case I(0) cancel each other out. As we mentioned, we can cut a composite region of the form I(1) in two different ways, so that the contributions from I(1) cancel each other out as well. We have illustrated some examples in Fig.~\ref{I1}.

By a simple exercise in planar geometry, we see that all the possible configurations from the contributions of I(2) and II(0) (up to cyclic permutations) can be illustrated as in Fig.~\ref{newII0} and they cancel each other out. Note that we draw the pictures in the simplest way; in general there can be some branch cuts and the regions would pass through the branch cuts and go to the other grid. Also in each picture the part of the pseudo-domain that has no corners in common with the empty rectangle, can be replaced with the other option. For example, in the left picture in Fig.~\ref{newII0}, $p$ is of the form L(2) and $p'$ is of the form L(3). Alternatively $p$ can be of the form L(1), i.e. the generator $\x$ can be such that $w_1\in \x$ and we do not need to use the L-shape ${L_1}{L_2}$. In this case the term $p*r$ from II(0) is again cancelled out against a term $r'*p'$ from I(2), but this time $p$ is trivial and $p'$ is of the form L(2), because we eliminated the L-shape ${L_1}{L_2}$ from $p$ and $p'$ in this alternative case.

Contributions from I(1$'$) will cancel out against contributions from I(3) and II(1); this can be seen by adding the appropriate column (either in the first or the second grid) that contains $O_1$ and $X_1$, to the composite region of I(1$'$). The configuration that we group with the contribution from I(1$'$) depends on the place of the component of the generator $\x$ on the lift of $\B_2$ in the column that we added. All the possible configurations are shown in Figs.~\ref{newI3L},\ref{newI3R} and \ref{newII1}. Here again, there is an alternative configuration for
each picture that we do not draw. Also there can be some branched cuts that the composite regions go through them, but we illustrated the configuration in the simplest way.

\end{proof}

Now we turn to prove that $F$ is a chain map with coefficients in $\Z$.

\begin{lemma} \label{FZ}
The map $F:C(\t H)\rightarrow C'$ preserves the Alexander grading, and is a chain map over $\Z$.
\end{lemma}

\begin{proof}
This lemma is a generalization of lemma~\ref{FZ2} and we only need to show that the configurations that are grouped together in the proof of lemma~\ref{FZ2}, have the right signs. By abuse of notation (for convenience), we show the formal rectangle $F(r)$ associated with an empty rectangle $r$, simply by $r$.

Depending on the order of the compositions in different cases, we need to prove the following equalities:
\begin{itemize}
\item If $p,p'\in\sigma^L$ and $p*r=r'*p'$, then $\mu(p)\cdot \SS(r)=\SS(r')\cdot \mu(p')$
\item If $p,p'\in\sigma^R$ and $p*r=r'*p'$, then $\mu(p)\cdot \SS(r)=-\SS(r')\cdot \mu(p')$
\item If $p,p'\in\sigma^F$ and $p*r=p'*r'$ then $\mu(p)\cdot \SS(r)=-\mu(p')\cdot \SS(r'),$ (similarly for $r*p=r'*p'$)
\end{itemize}
\noindent where $p,p'\in \sigma^F$, $r$ and $r'$ are empty rectangles.

Note that throughout this proof, for simplicity we call the rectangles and pseudo-domains with alphabets that shows their underlying regions. However, this notion does not specify the initial generators. In other words, we represent two rectangles with the same support and different initial generators, with the same alphabets. In this case, the initial generator of the rectangles should be understood from the context.
 
Note that $\SS$ depends on the initial formal generator of a rectangle or pseudo-domain. Let $A$ and $B$, be two formal rectangles, in our notation when we write $\SS(A)\cdot \SS(B)$, this means that the initial formal generator of $B$ is the same as the terminal formal generator of $A$, hence the terms do not freely commute.


\textbf{Contributions from I(0):}

In case of terms in I(0), the cancelling pairs are of the form $p*r=r'*p'$, where $r$ and $r'$ are of Type~1 and have the same support, also $p,p'\in \sigma^F$ have the same support but they differ in the initial generators; If $p,p'\in \sigma^L$ are of the form L(1), $\mu(p)=\mu(p')=1$. If $p,p'\in \sigma^L$ are of the form L(2) (resp. L(3)), we need to use the property S-3, two (resp. four) times, and this gives the desired equality.

Suppose that $p,p'\in \sigma^R$ are octagons. The sign of an octagon is defined by the product of the signs of the three associated formal rectangles. So we need to use the property S-3 of the sign assignments, three times. This gives a $(-1)^3=-1$.
If $p,p'\in \sigma^R$ are of the form R(2) (resp. R(3)), we need to use the property S-3, five (resp. seven) times, that gives a $-1$. Hence we get the desired equality.


\textbf{Contributions from I(1):}

Consider a rectangle $r$ and a pseudo-domain $p$, where they share a corner, and $r$ is disjoint from $w_0$ and $w_1$. Note that $p$ cannot be trivial in this case.

Let $p\in \sigma^F$ be a non-trivial pseudo-domain of Type L or R. Let $m$ be the number of formal rectangles that is used in defining the sign of $p$. Let $r_i$ ($i=1,\dots,m$) be the $i^{th}$ formal rectangle that we use for defining the sign of $p$.

\begin{figure}[h]
\centerline{\includegraphics[scale=0.7]{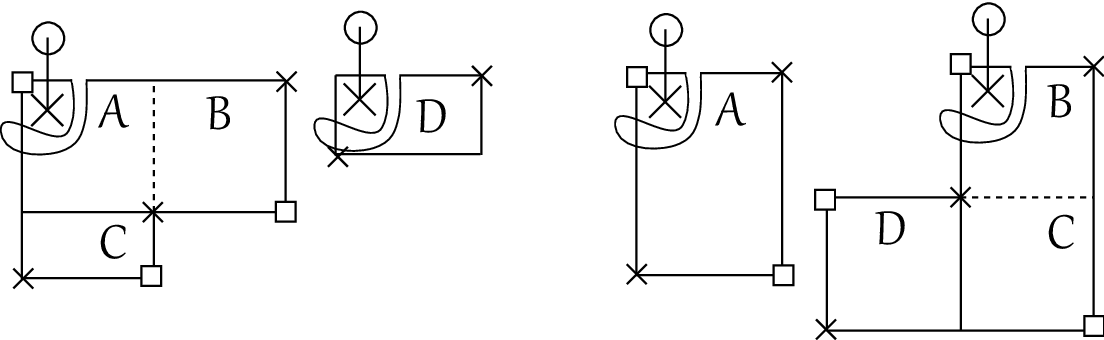}}
\caption {We have illustrated two examples of the pairs in I(1). The pseudo-domains $p,p'\in \sigma^R$, that we used are octagons. The left picture is of the form $r*p=r'*p'$ where $r=C$, $p=ABD$, $r'=B$ and $p'=ACD$. The right picture is of the form $r*p=p'*r'$ where $p=ABC$, $r=D$, $r'=CD$ and $p'=AB$.}
\label{octagonjuxta}
\end{figure}

\begin{itemize}

\item Suppose that the composite domain of $r$ and $p$ is of the form $r*p$ and the alternative decomposition is $r'*p'$. Among the pseudo-rectangles $r_i$ ($i=1,\dots,m$) that appear in the decomposition of $p$, one shares a corner with $r$, we call that $r_j$. 
If we start from the decomposition of $\SS(r)\cdot \mu(p)$, we need to use the property S-3, $j-1$ times to reach a decomposition such that a term of the form $\SS(s)\cdot \SS(r_j)$ appears, where $s$ has the same support as $r$. Then we need to use the property S-3, once more, to get $S(s')\cdot S(r'_j)$, where $s'$ has the same support as $r'$. At last, we need $j-1$ more use of the property S-3 to get to $\SS(r')\cdot \mu(p')$. Here we used the property S-3, $2(j-1)+1$ times. Hence we obtain $\SS(r)\cdot \mu(p) = - \SS(r')\cdot \mu(p')$. 
The case $p*r=p'*r'$ is similar. 

As an example, see the left picture of Fig.~\ref{octagonjuxta}. Here we illustrated a juxtaposition of a rectangle $C$ and an octagon $ABD$ that has an alternate decomposition $r*p=r'*p'$, where $r=C$, $p=ABD$, $r'=B$ and $p'=ACD$.
\begin{eqnarray*}
\SS(r)\cdot \mu(p)
&=&[\SS(C)\cdot \SS(r_{AB})]\cdot \SS(r_D)\cdot \SS(R_{ABD})\\
&=&-\SS(B)\cdot \SS(r_{AC})\cdot \SS(r_D)\cdot \SS(R_{ABD})\\
&=&-\SS(r')\cdot \mu(p').
\end{eqnarray*}
Note that the third formal rectangle in the definition of the sign of the octagons $ABD$ and $ACD$ are the same.

\item Suppose that the composition domain of $r$ and $p$ is of the form $r*p$, and the alternative decomposition is of the form $p'*r'$.
Among the pseudo-rectangles $r_i$ ($i=1,\dots,m$) that appear in the decomposition of $p$, one shares a corner with $r$, we call that $r_j$.  
If we start with the decomposition of $\SS(r)\cdot \mu(p)$, we need to use the property S-3, $j-1$ times, so that we reach a decomposition where a term $\SS(s)\cdot \SS(r_j)$ appears, here $s$ has the same support as $r$. Using the property S-3, once more, leads to a terms $\SS(r'_j)\cdot \SS(s')$, where $s'$ has the same support as $r'$. At last, by using the property S-3, $m-j$ times, we get to $(-1)^m\mu(p')\cdot \SS(r')$. If $p,p'\in \sigma^L$, $m$ is even. If $p,p'\in \sigma^R$, $m$ is odd. This gives the desired equality. The case $p*r=r'*p'$ is similar.

As an example, see the right picture of the Fig.~\ref{octagonjuxta}. The composite region has two decomposition of the forms $r*p=p'*r'$, where $p=ABC$, $r=D$, $r'=CD$ and $p'=AB$. In general, we denote the rectangles that we use to define the sign of an octagon $p=AB$, by $r_A$, $r_B$ and $R_{AB}$, respectively (i.e. $\mu(AB)=\SS(r_A)\cdot \SS(r_B)\cdot \SS(R_{AB})$).

Note that the third formal rectangle that we use to define the sign of the octagons $ABC$ and $AB$, have the same boundary arcs, the same orientations, and the same component of the generators in their corners. The only difference is the components of the initial generator which are not at the corners of the formal rectangles. We denote these formal rectangles by $R_{ABC}$ and $R_{AB}$, respectively.
\begin{eqnarray*}
\SS(r)\cdot \mu(p)
&=&\SS(D)\cdot \SS(r_A)\cdot \SS(r_{BC})\cdot \SS(R_{ABC})\\
&=&-\SS(r_A)\cdot \SS(D)\cdot \SS(r_{BC})\cdot \SS(R_{ABC})\\
&=&\SS(r_A)\cdot \SS(r_B)\cdot \SS(CD)\cdot \SS(R_{ABC})\\
&=&-\SS(r_A)\cdot \SS(r_B)\cdot \SS(R_{AB})\cdot \SS(CD)\\
&=&-\mu(p')\cdot \SS(r').
\end{eqnarray*}
In the fourth equality, the rectangle with the support $CD$ is disjoint from the formal rectangle $R_{ABC}$. When we switch their order, the initial generator of $R_{ABC}$ changes, and with the new initial generator, it becomes $R_{AB}$.

\end{itemize}

\textbf{Contributions of I(2) and II(0):}

Based on the argument in the proof of Lemma~\ref{FZ2}, we only need to obtain the sign assignment in the cases that we illustrated in Fig.~\ref{newII0}.

The left picture shows a term $p*r$ from II(0), where $p={L_1}{L_2}$ is a Type L pseudo-domain of the form L(2) and $r=ABCD$ is a Type 2 pseudo-rectangle, which cancels out with a term $r'*p'$ from I(2), where $r'=A$ is a Type 1 pseudo-rectangle and $p'=BCD{L_1}{L_2}$ is a Type L pseudo-domain of the form L(3).
\begin{eqnarray*}
\mu(p)\cdot \SS(r)
&=&\SS({L_1})\cdot \SS({L_2})\cdot \SS(A)\cdot \SS(BC)\cdot \SS(D)\\
&=&\SS(A)\cdot [\SS(BC)\cdot \SS(D)\cdot \SS({L_1})\cdot \SS({L_2})]\\
&=&\SS(r')\cdot \mu(p').
\end{eqnarray*}
\noindent In the second equality, since the terms with supports $L_1$ and $L_2$ are disjoint from the terms that we used to define the sign of $r$, we can move them to the right by using the property S-3 an even number of times.

The picture in the middle of the Fig.~\ref{newII0} shows a contribution $p*r$ from II(0) that cancels out with a contribution $p'*r'$ from I(2), where $p, p'\in \sigma^R$, $p=ABCDH$ is an octagon and $p'=ABDEFGH$ is of the form R(2), $r'=C$ is a Type 1 pseudo-rectangle and $r=AEFG$ is a Type 2 pseudo-rectangle. Note that the formal rectangle associated with $r$ has no edge in common with the three formal rectangles that we use for defining the sign of $p$. Hence in the second equality we move $r$, using the property S-3, three times. Note that the terms that we showed with $AEFG$, are actually two different rectangles with the same support. 
\begin{eqnarray*}
\mu(p)\cdot \SS(r)
&=&\mu(p)\cdot \SS(AEFG)\\
&=&-\SS(AEFG)\cdot \mu(ABCDH)\\
&=&-\SS(AEFG)\cdot \SS(r_{ABCD})\cdot \SS(r_H)\cdot \SS(R_{ABCDH})\\
&=&-\SS(AEFG)\cdot \SS(B)\cdot \SS(CD)\cdot \SS(r_{A})\cdot \SS(r_H)\cdot \SS(R_{ABCDH})\\
&=&\SS(AEFG)\cdot \SS(B)\cdot \SS(r_{AD})\cdot \SS(C)\cdot \SS(r_{H})\cdot \SS(R_{ABCDH})\\
&=&\SS(AEFG)\cdot \SS(B)\cdot \SS(r_{AD})\cdot \SS(r_{H})\cdot \SS(R_{ADH})\cdot \SS(C)\\
&=&\SS(AEFG)\cdot \SS(B)\cdot \mu(ADH)\cdot \SS(C)\\
&=&-\SS(ABG)\cdot \SS(EF)\cdot \mu(ADH)\cdot \SS(C)\\
&=&-\mu(p')\cdot \SS(r').
\end{eqnarray*}
\noindent In the fourth equality, we use Remark~\ref{rem} to write the sign of the formal rectangle $r_{ABCD}$ as $\SS(B)\cdot \SS(CD)\cdot \SS(r_A)$. In the sixth equality, when we move the term $\SS(C)$ using property S-3, the third formal rectangle $R_{ABCDH}$ in the definition of the sign of the octagon ${ABCDH}$, only differs with $R_{ADH}$ in the initial generators.

The right picture in Fig.~\ref{newII0}, shows a contribution $p*r$ from II(0) that cancels out with a contribution $r'*p'$ from I(2), where $p, p'\in \sigma^R$, $p=AGH$ is an octagon and $p'=ABDEFGH$ is of the form R(2), $r'=C$ is a Type 1 pseudo-rectangle and $r=ABCDEF$ is a Type 2 pseudo-rectangle.
\begin{eqnarray*}
\mu(p)\cdot \SS(r)
&=&\mu(AGH)\cdot \SS(ABCDEF)\\
&=&-\SS(ABCDEF)\cdot \mu(AGH)\\
&=&\SS(C)\cdot [\SS(ABF)\cdot \SS(DE)\cdot \mu(AGH)]\\
&=&\SS(r')\cdot \mu(p').
\end{eqnarray*}
\noindent In the second equality, the formal rectangle associated with $r$ has no edge in common with the three rectangles associated to the octagon $p$.
Hence we need to use the property S-3, three times.
Note that the two terms, with the support $ABCDDEF$, that are shown in the two sides of the second equality, are not the same rectangles, they differ in the initial generators. 

In the third equality, we can decompose the formal rectangle $ABCDEF$, using Remark~\ref{rem}, as follows:
$$\SS(ABCDEF)=\SS(C)\cdot \SS(ABF)\cdot \SS(DE)$$


\textbf{Contributions from I(1$'$) and I(3):}

We discussed in Lemma~\ref{FZ2}, that the contributions from I(1$'$) cancels out the contributions from I(3) and II(1). In the following we show the sign consistency in each of the cases in Figs.~\ref{newI3L},\ref{newI3R} and \ref{newII1}.

The left picture in Fig.~\ref{newI3L}, shows a contribution $r*p$ from I(3) that cancels out with a contribution $r'*p'$ from I($1'$), where $p, p'\in \sigma^L$, $p=ACD{L_1}{L_2}$ is of the form L(3) and $p'={L_1}{L_2}$ is of the form L(2), $r=B$ and $r'=D$ has $w_0$ as its lower-right corner.
\begin{eqnarray*}
\SS(r)\cdot \mu(p)
&=&[\SS(B)\cdot \SS(AC)]\cdot \SS(D)\cdot \SS(L_1)\cdot \SS(L_2)\\
&=&-\SS(D)\cdot \SS(L_1)\cdot \SS(L_2)\\
&=&-\SS(r')\cdot \mu(p').
\end{eqnarray*}
In the second equality we used the property S-2.

The right picture in Fig.~\ref{newI3L}, shows a contribution $p*r$ from I(3) that cancels out with a contribution $r'*p'$ from I($1'$), where $p, p'\in \sigma^L$, $p=ABE$ is of the form L(2) and $p'$ is of the form L(1), $r=CD$ is shaded horizontally and $r'=DE$ has $w_0$ as its lower-right corner.
\begin{eqnarray*}
\mu(p)\cdot \SS(r)
&=&\SS(AB)\cdot [\SS(E)\cdot \SS(CD)]\\
&=&-[\SS(AB)\cdot \SS(C)]\cdot \SS(DE)\\
&=&\SS(DE)\\
&=&\SS(r')\cdot \mu(p').
\end{eqnarray*}
Note that $p'$ is of Type L(1) and its sign is defined to be equal to $1$.


The left picture in Fig.~\ref{newI3R}, shows a contribution $r*p$ from I(3) that cancels out with a contribution $r'*p'$ from I($1'$), where $p, p'\in \sigma^R$, $p=ABCEFGH$ is of the form R(2) and $p'=ABH$ is an octagon, $r=D$ is shaded horizontally and $r'=FG$ has $w_0$ in the interior of its right edge.
\begin{eqnarray*}
\SS(r)\cdot \mu(p)
&=&[\SS(D)\cdot \SS(ACE)]\cdot \SS(FG)\cdot \mu(ABH)\\
&=&-\SS(FG)\cdot \mu(ABH)\\
&=&-\SS(r')\cdot \mu(p').
\end{eqnarray*}

The right picture in Fig.~\ref{newI3R}, shows a contribution $p*r$ from I(3) that cancels out with a contribution $r'*p'$ from I($1'$), where $p, p'\in \sigma^R$, $p=ABDFGHK$ is of the form R(2) and $p'=AHK$ is an octagon, $r=CE$ is shaded horizontally and $r'=EFG$ has $w_0$ in the interior of its right edge.
\begin{eqnarray*}
\mu(p)\cdot \SS(r)
&=&\SS(ABD)\cdot \SS(FG)\cdot [\mu(AHK)\cdot \SS(CE)]\\
&=&-\SS(ABD)\cdot [\SS(FG)\cdot \SS(CE)]\cdot \mu(AHK)\\
&=&[\SS(ABD)\cdot \SS(C)]\cdot \SS(EFG)\cdot \mu(AHK)\\
&=&-\SS(EFG)\cdot \mu(AHK)\\
&=&-\SS(r')\cdot \mu(p').
\end{eqnarray*}

In the second equality, the rectangle $CE$ is disjoint from the octagon $AHK$, hence we use the property S-3 three times.


\textbf{Contributions from I(1$'$) and II(2):}

The left picture in Fig.~\ref{newII1}, shows a term $p*r$ from II(1) that cancels out a term $r'*p'$ from I(1$'$), where $p, p' \in \sigma^R$, $p=ABF$ and $p'=AF$ are an octagons, $r=ACDE$ is a Type 2 pseudo-rectangle, and $r'=DE$ is shaded gray.
Fig.~\ref{octagonII1} is an illustration of the formal flows at hand. Note that the formal rectangles associated with $K$ and $DE$ are the same.
\begin{eqnarray*}
\mu(p)\cdot \SS(r)
&=&\SS(r_F)\cdot \SS(r_{AB})\cdot \SS(R_{ABF})\cdot \SS(ACDE)\\
&=&-\SS(r_F)\cdot [\SS(r_{AB})\cdot \SS(ACDE)]\cdot \SS(R_{AF})\\
&=&\SS(r_F)\cdot \SS(K)\cdot \SS(r_A)\cdot \SS(R_{AF})\\
&=&-\SS(K)\cdot \SS(r_F)\cdot \SS(r_A)\cdot \SS(R_{AF})\\
&=&-\SS(r')\cdot \mu(p').
\end{eqnarray*}
In the second equality, the formal rectangle associated with the Type~2 pseudo-rectangle $r$ has no edge in common with the formal rectangle $R_{ABF}$. By using the property S-3, we see that when we change the order of the terms, the terminal generator of the formal rectangle $ACDE$, which is the initial generator of the next term, is the same as the initial generator of the formal rectangle $R_{AF}$.
In the third equality, we use the property S-3.

\begin{figure}[h]
\centerline{\includegraphics[scale=0.8]{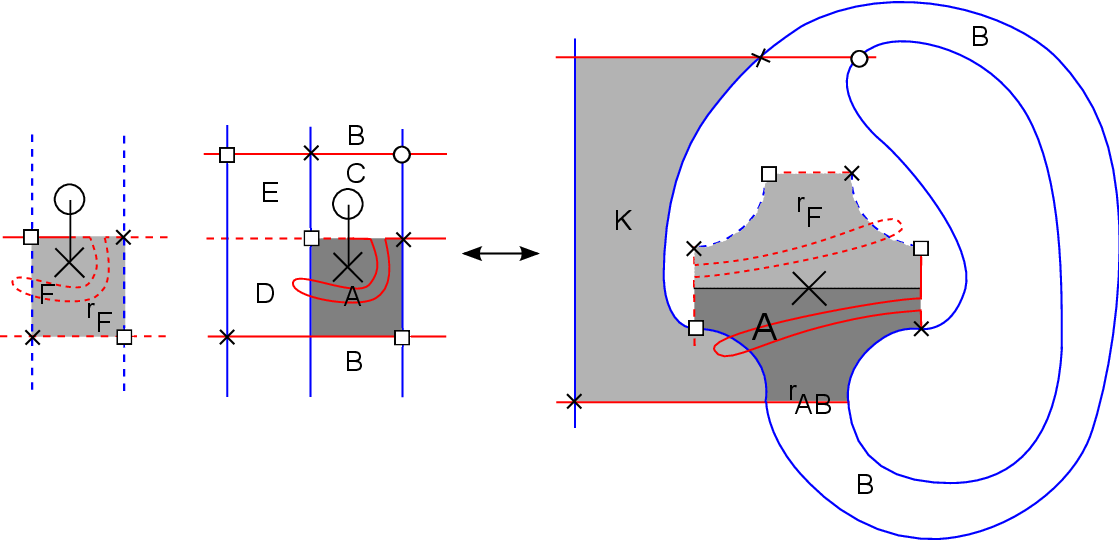}}
\caption {\textbf{A contributions from II(1) and I(1$'$)}
Here we illustrated the left picture from Fig.~\ref{newII1}. It shows a term $p*r$ from II(1) that cancels out a term $r'*p'$ from I(1$'$), where $p, p' \in \sigma^R$, $p=ABF$ and $p'=AF$ are an octagons, $r=ACDE$ is a Type 2 pseudo-rectangle, and $r'=DE$ is a Type 1 pseudo-rectangle.}
\label{octagonII1}
\end{figure}

The picture in the middle of the Fig.~\ref{newII1}, shows a term $p*r$ from II(1) that cancels out a term $r'*p'$ from I(1$'$), where $p, p' \in \sigma^R$, $p=ACDEGHJKLM$ is of the form R(2), $p'=CDEGHJKLM$ is of the form R(2) and is in gray, $r=BCDEFGHI$ is a Type 2 pseudo-rectangle, and $r'=FGHI$ is shaded horizontally.
\begin{eqnarray*}
\mu(p)\cdot \SS(r)
&=&\SS(ACDE)\cdot \SS(GHJK)\cdot [\mu(DLM)\cdot \SS(BCDEFGHI)]\\
&=&-\SS(ACDE)\cdot [\SS(GHJK)\cdot \SS(BCDEFGHI)]\cdot \mu(DLM)\\
&=&\SS(ACDE)\cdot \SS(BCDEFGHI)\cdot \SS(GHJK)\cdot \mu(DLM)\\
&=&\SS(ACDE)\cdot \SS(B)\cdot \SS(FGHI)\cdot \SS(CDE)\cdot \SS(GHJK)\cdot \mu(DLM)\\
&=&-\SS(FGHI)\cdot [\SS(CDE)\cdot \SS(GHJK)\cdot \mu(DLM)]\\
&=&-\SS(r')\cdot \mu(p').
\end{eqnarray*}
Note that the formal rectangle associated with the Type~2 pseudo-rectangle $r$ has no edge in common with the formal rectangles associated with the octagon $DLM$. Hence using the property S-3 three times, we obtain the second equality.
In the fourth equality, we use the Remark~\ref{rem} to write the sign of a formal rectangle $BCDEFGHI$, whose associated embedded rectangle in the Heegaard diagram has a component of the initial generator in its interior, as follows:
$$\SS(BCDEFGHI)=\SS(B)\cdot \SS(FGHI)\cdot \SS(CDE).$$
Note that the definition of the sign depends on the place of the component of the initial generator in the interior of the rectangle with support $BCDEFGHI$. Since the term with support $ACDE$ comes first, the rectangles in the decomposition of $r$ meet in $v_0$.

The right picture in Fig.~\ref{newII1}, shows a term $p*r$ from II(1) that cancels out a term $r'*p'$ from I(1$'$), where $p, p' \in \sigma^L$ are of the form L(3), $p=ACEG{L_1}{L_2}$, $p'=CEG{L_1}{L_2}$, $r=BCDEF$ is a Type 2 pseudo-rectangle, and $r'=DEF$ is a Type 1 pseudo-rectangle that has $w_1$ in the interior of its right edge.
\begin{eqnarray*}
\mu(p)\cdot \SS(r)
&=&\SS(L_1)\cdot \SS(L_2)\cdot \SS(AC)\cdot [\SS(EG)\cdot \SS(F)]\cdot \SS(BC)\cdot \SS(DE)\\
&=&-\SS(L_1)\cdot \SS(L_2)\cdot \SS(AC)\cdot \SS(EF)\cdot [\SS(G)\cdot \SS(BC)]\cdot \SS(DE)\\
&=&\SS(L_1)\cdot \SS(L_2)\cdot \SS(AC)\cdot \SS(EF)\cdot \SS(BC)\cdot [\SS(G)\cdot \SS(DE)]\\
&=&-\SS(L_1)\cdot \SS(L_2)\cdot \SS(AC)\cdot \SS(EF)\cdot [\SS(BC)\cdot \SS(D)]\cdot \SS(EG)\\
&=&\SS(L_1)\cdot \SS(L_2)\cdot \SS(AC)\cdot \SS(EF)\cdot \SS(BD)\cdot \SS(C)\cdot \SS(EG)\\
&=&\SS(AC)\cdot [\SS(EF)\cdot \SS(BD)]\cdot [\SS(L_1)\cdot \SS(L_2)\cdot \SS(C)\cdot \SS(EG)]\\
&=&-[\SS(AC)\cdot \SS(B)]\cdot \SS(DEF)\cdot \mu(p')\\
&=&\SS(r')\cdot \mu(p').
\end{eqnarray*}
Note that using the property S-2 in the last equality, we have $\SS(AC)\cdot \SS(B)=-1$, since the formal rectangle associated with the punctured rectangle $AC$, makes a Type $\B$ degeneration together with the formal rectangle associated with $B$.

\end{proof}


We decompose the set of generators of the Heegaard diagram $\t H$ according to the position of the components of a generator on $\t\alpha^0_1$ and $\t\alpha^1_1$. 
We represent the type of each generator with an ordered pair. The first (resp. second) entry is $I$ if the component on $\t\alpha^0_1$ (resp. $\t\alpha^1_1$) is on one of the lifts of $\B_1$ (either $\t\beta_1^0$ or $\t\beta_1^1$). The first (resp. second) entry is $J$ when the $\t\alpha_1^0$ (resp. $\t\alpha_1^1$) component is on one of the lifts of $\beta_2$. The $N$ in the first (resp. second) entry shows that the generator has its $\t\alpha^0_1$ (resp. $\t\alpha^1_1$) component neither on the lifts of $\beta_1$ nor on the lifts of $\beta_2$. Hence we have the following decomposition for the set of generators:
\[\S(\t H)=(I,I)\cup(I,J)\cup(I,N)\cup(J,I)\cup(N,I)\cup(J,J)\cup(J,N)\cup(N,J)\cup(N,N)\]

Having the above decomposition of the set of generators of $C$ we get a decomposition of $C$ as the direct sum of the sub-modules generated by the associated generators. 
\[C=C^{I,I}\oplus C^{I,J}\oplus C^{I,N}\oplus C^{J,I}\oplus C^{N,I}\oplus C^{J,J}\oplus C^{J,N}\oplus C^{N,J}\oplus C^{N,N}\\ \]

In order to show that $F$ is a quasi-isomorphism, following the ideas of \cite{MOST}, we define a filtration.

Denote by $\S(\t H,k)$ the set of generators $\x\in S(\t H)$ with Alexander gradings equal to $k\in \Z$. Let $C(\t H,k)$ be the complex generated by elements of $\S(\t H,k)$. Since $\partial$ preserves the Alexander grading $C(\t H,k)$ is a summand of $C(\t H)$.

Let $\Pi:\t H\rightarrow H$ be the double branched cover map and $Q\subset H$ be a set consisting of one point in each square of $H$ other than those in the same row or the same column as $\Pi(O_1)$.
Let $\x,\y\in \S(\t H)$ and $p\in \sigma(\x,\y)$, the image of $p$ under $\Pi$ is a 2-chain (not necessarily a domain) in the Heegaard diagram $H$. We denote by $X_1(p)$ (respectively $O_1(p)$), the number of times that $\Pi(X_1)$ (respectively $\Pi(O_1)$) lies in $\Pi(p)$ counted with multiplicity.

\begin{lemma}\label{2chain}
Let $D$ be a 2-chain (not necessarily a domain) in the Heegaard diagram $H$ that contains no basepoints, such that $\partial D$ consists of $\A$- and $\B$-circles. Then $D$ is trivial. ($D$ might contain a basepoint in its interior.)
\end{lemma}

\begin{proof}
Let $D$ be as above. We add some number of rows and some number of columns to $D$ to get a 2-chain $D'$, so that each square has non-negative multiplicity. There might be some rows or columns in $D'$ i.e. a row or column in $H$ with positive multiplicities in $D'$, we subtract them to get a 2-chain $D''$, with non-negative multiplicities and containing no row or column.

If we consider two adjacent squares in $H$, if they have multiplicities $a$ and $b$ in $D''$, then the segment between them will occur in the boundary of $D''$ with multiplicity $a-b$ (Note that the sign of $a-b$ reflects the direction in which this segment appears in the boundary of $D''$, and if $a=b$ then this segment has multiplicity $0$ in the boundary). Hence if $\A_1$ is a circle in the boundary of $D''$, since all the squares have non-negative multiplicity in $D''$, then depending on the orientation of $\A_1$ that comes from $D''$, one of the rows that have $\A$ as their boundary, has positive multiplicity in $D''$ which contradicts the way that we defined $D''$. More precisely if the induced orientation of $\A_1$ from the boundary of $D''$ is from left to right (or right to left), for any segment of it if the multiplicities of the squares above and below it are $a$ and $b$ then we will have $a \geq a-b=mult(\alpha_1) > 0$  (or $b \geq b-a=mult(\alpha_1) > 0$), which holds for the multiplicities of all the squares on the row above (or below) $\A_1$ depending on the orientation of $\A_1$. So there is no $\A$-circle in $\partial D''$. Similarly we can prove that there is no $\B$-circle in the boundary of $D''$. One can see each pair of squares that has an edge in common have the same multiplicity in $D''$.
From the definition of $D''$ we know that there exist squares that have zero multiplicity in $D''$, since it does not contain any row or column. So $D''$ is trivial. Hence $D$ is a sum of a number of rows and a number of columns.

In the representation of this two-chain as sum of rows and columns if a row has multiplicity $c$, then it contains an $O_i$ for some $i$. Since this two-chain does not contain $O_i$, the column through $O_i$ should have multiplicity $-c$. Now this column contains an $X_j$ so the row through that should have multiplicity $c$. Continuing in this way, since we work with a knot rather than a link, we see that all the rows have multiplicity $c$ and all the columns have multiplicity $-c$. So the two-chain is in fact trivial.

\end{proof}

For fixed $\x , \y \in \S(\t H,k)$, if $p, p' \in \sigma (\x , \y )$ are two pseudo-domains such that $O_1(p)=X_1(p)=O_1(p')=X_1(p')$ and there are no other basepoints inside them, counting with multiplicity by Lemma~\ref{2chain}, we have: $$\# (Q\cap \Pi(p))=\# (Q\cap \Pi(p')).$$

We call a pseudo-domains a \emph{Q-fine} pseudo-domains if its multiplicity at each base point except (possibly) $X_1$ is zero. So we can find a function $\QQ$ such that
\begin{equation}\label{QQ}
\QQ(\x)-\QQ(\y)=\# (Q\cap \Pi(p))
\end{equation}
\noindent for all Q-fine pseudo-domains $p \in \sigma (\x, \y)$.

The construction of $\QQ$ is a simple combinatorial process. Fix a generator $\x_{0}\in \S(\t H)$ and assign an arbitrary value to $\QQ(\x_0)$. Let $\y\in \S(\t H)$, we can define the value of $\QQ$ at $\y$ if there is a Q-fine pseudo-domain between $\x_0$ and $\y$ (i.e. from $\x_0$ to $\y$ or from $\y$ to $\x_0$.) If there is a Q-fine pseudo-domain $p\in \sigma(\x_0,\y)$, then according to Equation~\ref{QQ} the value of $\QQ$ on $\y$ is defined to be $\QQ(\y) = \QQ(\x_0)-\# (Q\cap \Pi(p))$. If there is a Q-fine pseudo-domain $p'\in\sigma (\y , \x_0 )$, we have $\QQ(\y)=\QQ(\x_0)+ \# (Q\cap \Pi(p))$.
Similarly if there is a chain of Q-fine pseudo-domains between $\x_0$ and $\y$, i.e. there are a number of generators $\x_i$ for $i=1,\cdots, n$, where $\x_n=\y$, and there is a Q-fine domain between each two consecutive generators $\x_i$ and $\x_{i+1}$. Then the value of $\QQ$ on $\y$ can be defined. If there is no such chain between a generator $\x'$ and $\x_0$ we can define $\QQ(\x')$ arbitrarily and continue the above procedure. Since the number of generators is finite this procedure will stop at some point.

In order to show that $\QQ$ is well-defined we have to consider two such chains between $\x_0$ and $\y$, and show that using each chain we get the same value for $\QQ$ at $\y$. To put it differently if we consider the union of these chains, we get a loop starting and ending at $\x_0$ and we want to show that if we follow the path the changes in the value of $\QQ$ adds up to zero. But for each edge of this loop, the change in the value of $\QQ$ along this edge is equal to $\# (Q\cap \Pi(p_i))$, where $p_i$ is the pseudo-domain associated to that edge. Hence the sum of changes in the value of $\QQ$ along this path is equal to $\# (Q\cap \sum(-1)^{n_i}\Pi(p_i))$ where $n_i$ is 1 or 0 depending on whether $p_i$ is in $\sigma (\x_i , \x_{i+1} )$ or in $\sigma (\x_{i+1} , \x_{i})$. 

Note that $p=\sum(-1)^{n_i} p_i$ is a pseudo-domain from a generator to itself so its boundary consists of a number of copies of $\alpha$- and $\beta$-circles. Note that since each $p_i$ has no basepoints except possibly $X_1$, the same holds for their sum that is $p$. Now $p$ is from a generator to itself so it does not change the Alexander grading, so the multiplicity of $X_1$ should also be zero. Hence when we take the image of $p$ under the $\Pi$ we get a two-chain satisfying the condition of the Lamma~\ref{2chain}, hence $\Pi(p)$ is trivial and in particular this shows that $\# (Q\cap p)= \# (Q\cap\sum(-1)^{n_i}\Pi(p_i))=0$. Hence the total changes of $\QQ$ along this loop is zero, which shows that $\QQ$ is well-defined.  One should note that in this proof we see that the two-chain obtained by the projection $\Pi$ is trivial, which does not mean that the domain in the Heegaard diagram of $\t H$ is trivial.

Using $\QQ$ we can define a filtration on $C(\t H,k)$. The boundary operator of the associated graded object counts those rectangles 
which do not contain any basepoints and does not contain any points of $Q$. We denote by $C_Q$ the associated graded object, and omit the index $k$ when there is no confusion.

In order to study the boundary map between different submodules and calculate the homology groups, we need a new definition. 

\begin{definition}
Given $\x, \y \in S(\t H)$, we call an empty rectangle $r \in Rect(\x,\y)$ that has no points of $Q$ inside it, an \emph{undone} rectangle with respect to the initial generator $\x$. In this case, the components of $\x$ are the upper-right and lower-left corners of $r$. We call a rectangle $r \in Rect(\y,\x)$ that has no points of $Q$ inside it, a \emph{done} rectangle with respect to $\x$.
See Fig.~\ref{cNN}.

Each term in the $\partial\x$ is obtained by picking an undone rectangle $r_0 \in Rect(\x,\y)$. We refer to the generator $\y$, as \emph{the generator obtained from  $\x$ by using $r_0$}.
\end{definition}

\begin{lemma} \label{Homology}
$H_*(C_Q)$ is isomorphic to free $\Z_2$-module generated by elements of $(I,I)$ and $(J,J)$. 
\end{lemma}  

\begin{proof}
There are two cases depending on whether $X_2$ is placed in the square, just to the left or just to the right of $O_1$, in the grid diagram $H$ for $K$. First, suppose $X_2$ is in the square just to the left of the square marked $O_1$.
Then we have a direct splitting 
$\Ctq=\Ctq^{I,I} \oplus B_1 \oplus B_2 \oplus B_3$, where $B_1$ is the module generated by all generators of types $(I,N)$ and $(I,J)$, $B_2$ is the module generated by all generators of types $(N,I)$ and $(J,I)$, $B_3$ is the module generated by all generators of types $(N,N)$, $(J,N)$, $(N,J)$ and $(J,J)$, $B_4$ is the module generated by all generators of types $(N,N)$, $(J,N)$ and $(N,J)$.

The differentials in $\Ctq^{I,I}$ are trivial, hence its homology is the free 
$\Z_2$-module generated by elements of $(I,I)$. Also $B_1$, $B_2$ and $B_3$ are chain complexes fitting into these exact sequences 

\begin{center}
\exact{\Ctq^{I,N}}{B_1}{\Ctq^{I,J}}

\exact{\Ctq^{N,I}}{B_2}{\Ctq^{J,I}}

\exact{B_4}{B_3}{\Ctq^{J,J}}\label{seq3}
\end{center}

Here $B_4$ fits into the below short exact sequence
\begin{center}
\exact{\Ctq^{N,N}}{B_4}{\Ctq^{J,N}\oplus \Ctq^{N,J}}
\end{center}

For explanation of these sequences see the captions of Figs.~\ref{cNN}, \ref{cIJ} and \ref{cJN}.
First we show that $H_*(\Ctq^{N,N})$ is zero. 
Note that for $\x \in S(\t H)$ we can consider all the done and undone rectangles, then the terms of $\partial \x$ in $\Ctq$ are given by counting undone rectangles with respect to $\x$. 

In the first Heegaard diagram in Fig.~\ref{cNN}, there are exactly three disjoint undone rectangles and one done rectangle, with respect to $\x$. We represent $\x$ by its undone rectangles, $e_1\wedge e_{1}' \wedge e_2$.
In each term in $\partial \x$ in $\Ctq$, one of these undone rectangles (with respect to $\x$) will be counted. For example, we denote by $ e_1\wedge e_{1}'$, the term that comes from using $e_2$. In other words, it represents a generator $\y$ such that $e_1$ and $e_{1}'$ are the only undone rectangles with respect to $\y$. Hence $\partial \x$ in $\Ctq$ can be represented by the following notation (we call this the induced boundary operator):
\[
\partial (e_1\wedge e_{1}' \wedge e_2) = e_1\wedge e_{1}' + e_{1}' \wedge e_2 + e_1 \wedge e_2.
\]

One can see that each generator in $\Ctq$ can have at most four undone rectangles. Given a generator $\x$, we construct a set of generators associated with $\x$. This set consists of the generators obtained from $\x$ using a number of undone rectangles, and also those generators that $\x$ is obtained from them by using a number of undone rectangles. By the above notation, this set with the induced boundary operator from $C_Q$ is isomorphic to an exterior algebra over a vector space of dimension at most four with coefficients in $\Z_2$.

In the above case, the vector space is four dimensional.
As another example, see the second Heegaard diagram in Fig.~\ref{cNN}; there are two disjoint undone rectangles with respect to $\x$. Hence we represent the generator $\x$ with $e_1 \wedge e_2$. The induced boundary map is represented as follows: 
\[
\partial (e_1 \wedge e_2) = e_1 + e_2
\]
Note that here the exterior algebra is over a three dimensional vector space.

This shows that we can partition our complex as a union of subcomplexes, each isomorphic to the exterior algebra over a vector space with coefficients in $\Z_2$. Since this complex is isomorphic to the complex of the reduced homology of an $n$-simplex over $\Z_2 $, hence the homology of $\Ctq^{N,N}$ being equal to the direct sum of the homology of these subcomplexes is trivial. The same idea shows that $H_*(\Ctq^{I,N})=0$ and $H_*(\Ctq^{N,I})=0$.

\begin{figure}[h]
\centerline{\includegraphics[scale=0.8]{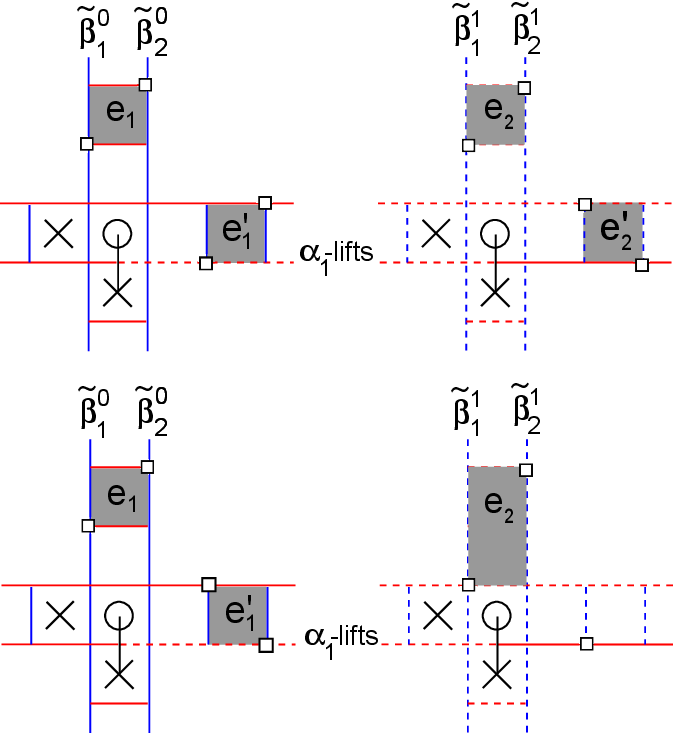}}
\caption {Hollow squares show a generator $\x \in \Ctq^{N,N}$. In the first Heegaard diagram, the rectangle $e_{2}'$ is done with respect to $\x$, and the rectangles $e_1, e_{1}', e_2 $ are undone with respect to $\x$. In the second Heegaard diagram $e_1'$ is done with respect to $\x$, and $e_1$, $e_2$ are undone rectangles with respect to $\x$. Note that $\A_1^0$ is drawn with a solid line, and $\A_1^1$ is drawn with a dashed line.}
\label{cNN}
\end{figure}

\begin{figure}[h]
\centerline{\includegraphics[scale=0.8]{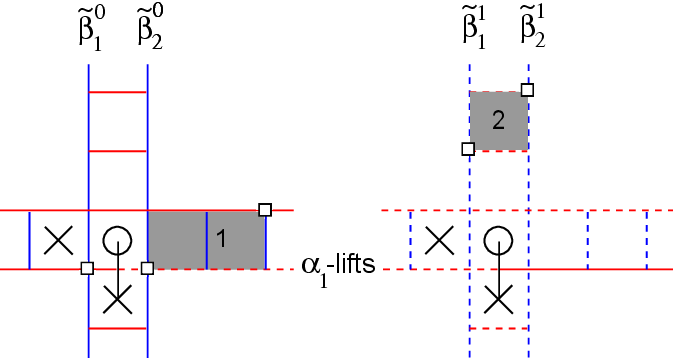}}
\caption {Hollow squares show a generator $\x \in \Ctq^{I,J}$. We count the rectangles number 1 and 2, in the boundary map. Using rectangle 1 leads to a generator in $\Ctq^{I,N}$. If we use rectangle 2 the result is a generator in $\Ctq^{I,J}$. Here the circle $\A_1^0$ is drawn in solid, and the circle $\A_1^1$ is drawn with a dashed line.}
\label{cIJ}
\end{figure}

In order to compute the homology of 
$B_1$, it is enough to compute the homology of $\Ctq^{I,J}$ as a quotient.
Each generator in $\Ctq^{I,J}$ has at most two terms in its boundary, one of which is in $\Ctq^{I,N}$. Hence $H_*(B_1)=H_*(\Ctq^{I,J})=0$, See Fig.~\ref{cIJ}. In the same way we can prove that $H_*(B_2)=H_*(\Ctq^{J,I})=0$ and 
$H_*(B_4)=H_*(\Ctq^{J,N}\oplus \Ctq^{N,J})=0$, See Fig.~\ref{cJN}. Finally the differentials in $\Ctq^{J,J}$ are trivial, so its homology is the 
free $\Z_2$-module generated by the elements of $(J,J)$.

\begin{figure}[h]
\centerline{\includegraphics[scale=0.8]{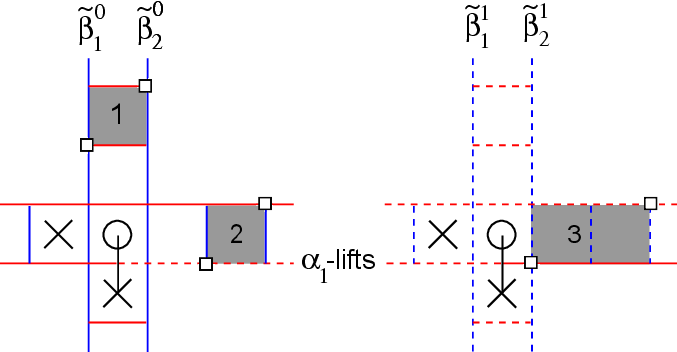}}
\caption {Hollow squares show a generator $\x \in \Ctq^{J,N}$. We count each of the rectangles shown above, in the boundary map. Counting rectangles 1 or 2 lead to a generator in $\Ctq^{J,N}$. If we use rectangle 3, the result is a generator in $\Ctq^{N,N}$. Here $\A_1^0$ (resp. $\A_1^1$) is drawn in solid (resp. dashed) line.}
\label{cJN}
\end{figure}

Suppose on the other hand that $X_2$ is just to the right of $O_1$, in the grid diagram $H$ for $K$. Then there is a direct sum splitting 
$\Ctq=\Ctq^{J,J} \oplus B'_1 \oplus B'_2 \oplus B'_3$, where $B'_1$ is the module generated by all generators of types $(J,I)$ and $(J,N)$, $B'_2$ is the module generated by all generators of types $(I,J)$ and $(N,J)$, $B'_3$ is the module generated by all generators of types $(I,I)$, $(I,N)$, $(N,I)$ and $(N,N)$, $B'_4$ is the module generated by all generators of types $(I,I)$, $(I,N)$ and $(N,I)$. Also $B'_1$, $B'_2$ and $B'_3$ are chain complexes fitting into these exact sequences 

\begin{center}
\exact{\Ctq^{J,I}}{B'_1}{\Ctq^{J,N}}

\exact{\Ctq^{I,J}}{B'_2}{\Ctq^{N,J}}

\exact{B'_4}{B'_3}{\Ctq^{N,N}} 
\label{2seq2}
\end{center}
Here $B'_4$ fits in the following exact sequence 
\begin{center}
\exact{\Ctq^{I,I}}{B'_4}{\Ctq^{I,N}\oplus\Ctq^{N,I}} 
\label{2seq3} 
\end{center}
All the differentials in $\Ctq^{J,J}$ are trivial, so its homology is the $\Z_2$-module generated by the elements of $(J,J)$.
It is easy to see that $H_*(\Ctq^{J,I})=0$. If we consider $\Ctq^{J,N}$ as a quotient chain complex, it is also easy to show that 
$H_*(\Ctq^{J,N})=0$.
Hence we have $H_*(B'_1)=0$.
In the same way we can prove that $H_*(B'_2)=0$.
The differentials in $\Ctq^{I,I}$ are trivial, so its homology is the $\Z_2$-module generated by the elements of $(I,I)$.
Also $H_*(\Ctq^{I,N}\oplus\Ctq^{N,I})=0$ as a quotient chain complex.
A simple calculation shows that $H_*(\Ctq^{N,N})=0$. Considering the above exact sequences, 
we get the desired result.
\end{proof}

\begin{proposition} \label{quasi}
The map $F$ is a filtered quasi-isomorphism.
\end{proposition}

Note that in order to define the function $\QQ$ for a generator $\x$ of $C(\t G)$, we consider $\psi(\x)$ and define $\QQ(\x):= \QQ(\psi(\x))$. For an element $(\x_0,\x_1) \in C'$ we define $$\QQ(\x_0,\x_1):=max\{\QQ(\x_0),\QQ(\x_1) \}.$$
 
\begin{proof}
First we show that $F$ preserve the $Q$-filtration. Let $\x$ be a generator of $C$ we want to show $\QQ(\x) \geq \QQ(F(\x))$. So we have to show that $\QQ(\x) \geq \QQ(F^L(\x))$ and $\QQ(\x) \geq \QQ(F^R(\x))$. 

We decompose a pseudo-domain of type $L$ or $R$ as the $*$ of a number of punctured rectangles, octagons that contain $X_1$ and empty rectangles.

Let $\mathbf{u}, \mathbf{v}\in \S(\t H)$ such that there is a punctured rectangle $\a \in A(\mathbf{u} , \mathbf{v})$, then the image of the complementary rectangle $r_{\a}\in Rect(\mathbf{v}, \mathbf{u})$ under $\Pi$ is empty from the points in $Q$. Hence by Equation~\ref{QQ} we have $\QQ(\mathbf{u})=\QQ(\mathbf{v})$.
For $\mathbf{w}, \mathbf{z}\in \S(\t H)$ if there is an octagon that contains $X_1$ or an empty rectangle from $\mathbf{w}$ to $\mathbf{z}$, from Equation~\ref{QQ} we have $\QQ(\bm w) \geq \QQ(\bm z)$. Hence $\QQ(\x) \geq \QQ(F(\x))$ i.e. $F$ preserves the $Q$-filtration.
\end{proof}

We consider the map induced by $F$:
\[F_Q:\Ctq \longrightarrow\Ctq '\]

By Lemma~\ref{Homology}, the homology of $\Ctq$ is carried by the subcomplex $\Ctq^{(I,I)}\oplus\Ctq^{(J,J)}\subset\Ctq$. So we consider the restriction of $F_Q$ to this subcomplex, and will show that it induces an isomorphism. Also by definition we have $\Ctq '=\L_Q\oplus\RR_Q$. Note that the only pseudo-domains that are allowed, are those that do not have intersection with the dots in $Q$.

Let $\x$ be a generator in $\Ctq^{(I,I)}$. Since $F(\x)$ ends up in $(I,I)$, the only pseudo-domain that contributes to $F(\x)$ is the trivial domain of Type $L$. Hence $F_Q^L$ restricted to $C_Q^{(I,I)}$ is an isomorphism.

The restriction of $F_Q^R$ to $\Ctq^{(J,J)}$, counts pseudo-domains of Type $R$ which are supported in the two rows and the two columns through $O_1$; that is we only count octagons.
But for each element in $(J,J)$ there is a unique way of assigning an octagon. Thus, in this case, $F_Q$ is a quasi-isomorphism.

Now we need the following algebraic lemma; for a proof see Theorem 3.2 from \cite{Mc}.

\begin{lemma} \label{nice}
Suppose that $F:C\longrightarrow C'$ is a filtered chain map which induces an isomorphism on the homology of the associated graded object. Then $F$ is a quasi-isomorphism. 
\end{lemma}

We now use Lemma~\ref{nice}, first for the filtration of $Q$, and then for the Alexander grading to conclude that $F$ is a quasi-isomorphism.

This completes the proof of the fact that $H(C)\cong H(C')\cong H(B)\otimes V$, where $V\cong \Z_2 \oplus \Z_2$ with generators in gradings $0$ and $-1$.


\begin{lemma}\label{gauge}
$H_*(C_Q)$ is isomorphic to free $\Z$-module generated by elements of $(I,I)$ and $(J,J)$. 
\end{lemma}  

\begin{proof}
This lemma is the generalization of Lemma~\ref{Homology} over $\Z$.
Here we use the same exterior algebra notation. For example if we start with generator $\x$ as in Fig.~\ref{cNN}, which is denoted by $e_1\wedge e_2 \wedge e_3$ and use the empty rectangle $e_3$, then we get the term $e_1 \wedge e_2$ with the sign $s_3$ of the formal rectangle associated to $e_3$ with the initial generator $\x$. So we get an expression of the form:
$$\partial  (e_1\wedge e_2 \wedge e_3) = s_1\cdot e_2 \wedge e_3 + s_2\cdot e_1 \wedge e_3 +s_3\cdot e_1\wedge e_2 $$
It is not hard to see that since the sign assignment has the property S-3, with an appropriate gauge transform, we can change the sign assignment such that we have:
$$\partial (e_{n_1}\wedge \cdots \wedge e_{n_k})= \displaystyle \sum (-1)^{i}  (e_{n_1}\wedge \cdots \wedge \widehat{e_{n_i}}\wedge \cdots\wedge e_{n_k})$$
This means that similar to case of $\Z_2$, we split the complex as the direct sum of a number of subcomplexes, each isomorphic to an exterior algebra of a vector space of dimension at most $4$ over $\Z$, and similarly they are isomorphic to the complex of reduced homology for an $n$-simplex ($n \leq 4$) over $\Z$. Hence the homology of the complex is the direct sum of the homology of these subcomplexes, and so it is trivial. 
\end{proof} 

Using Lemma~\ref{nice} and a similar argument as in Proposition \ref{quasi}, show that $H(C)\cong H(C')$ over $\Z$ and this concludes the proof of the fact that $H(C)\cong H(B) \otimes V$, where $V \cong \Z \oplus \Z$ with generators in gradings $0$ and $-1$.


\section*{Acknowledgments}
I am grateful to Zolt\'an Szab\'o for suggesting this problem, numerous helpful discussions, continuous advice through the course of this work, and reading a draft of this paper. I would also like to thank Iman Setayesh for helpful conversations and his contribution in defining the $\QQ$ filtration in Section~\ref{stab} and Lemma~\ref{gauge}. This work was done when I was a visiting student research collaborator at Princeton university, and I am grateful for the opportunity. This paper is a part of my Ph.D. thesis as a graduate student in Sharif University of Technology under the supervision of Mohammadreza Razvan. I also thank the referee for the helpful comments and suggestions, and Eaman Eftekhary for reading a draft of this paper.


\end{document}